\documentclass[12pt,letter]{amsart}
\usepackage{amsmath,amssymb}
\title[Factors
with at most one Cartan subalgebra]{On a class
of $\mathrm{II}_1$ factors\\
with at most one Cartan subalgebra}
\author[N. Ozawa]{Narutaka OZAWA$^*$}
\address{Department of Mathematical Sciences,
University of Tokyo, Komaba, \mbox{153-8914}\\
\indent Department of Mathematics, UCLA, Los Angeles, CA90095}
\email{narutaka@ms.u-tokyo.ac.jp}
\thanks{${}^*$ Supported by NSF-Grant and Sloan Foundation}
\author[S. Popa]{Sorin Popa$^{**}$}
\address{Department of Mathematics,
UCLA, Los Angeles, CA90095}
\email{popa@math.ucla.edu}
\thanks{${}^{**}$ Supported by NSF-Grant}
\dedicatory{Dedicated to Alain Connes on his 60th birthday.}

\date{June 27, 2007}
\newtheorem{thm}{Theorem}[section]
\newtheorem{lem}[thm]{Lemma}
\newtheorem{prop}[thm]{Proposition}
\newtheorem{cor}[thm]{Corollary}
\newtheorem*{thm*}{Theorem}
\newtheorem{cor*}{Corollary}
\theoremstyle{definition}
\newtheorem{defn}[thm]{Definition}
\newtheorem{rem}[thm]{Remark}
\numberwithin{equation}{section}
\newcommand{\R}{{\mathbb R}}
\newcommand{\C}{{\mathbb C}}
\newcommand{\N}{{\mathbb N}}
\newcommand{\Z}{{\mathbb Z}}

\newcommand{\F}{{\mathbb F}}
\newcommand{\B}{{\mathbb B}}
\newcommand{\GG}{{\mathcal G}}
\newcommand{\Sg}{{\mathcal S}}
\newcommand{\NN}{{\mathcal N}}
\newcommand{\Nor}{{\mathcal N}}
\newcommand{\Zt}{{\mathcal Z}}
\newcommand{\U}{{\mathcal U}}
\newcommand{\G}{\Gamma}
\newcommand{\e}{\varepsilon}
\newcommand{\op}{\mathrm{op}}
\newcommand{\p}{\varphi}
\newcommand{\hh}{{\mathcal H}}

\newcommand{\cb}{\mathrm{cb}}
\DeclareMathOperator{\Ad}{Ad}
\DeclareMathOperator{\Aut}{Aut}

\newcommand{\vt}{\mathbin{\bar{\otimes}}}

\DeclareMathOperator*{\Lim}{Lim}
\newcommand{\id}{\mathrm{id}}
\DeclareMathOperator{\Tr}{Tr}
\newcommand{\SL}{\mathrm{SL}}
\newcommand{\ip}[1]{\mathopen{\langle}#1\mathclose{\rangle}}
\newcommand{\oneplus}{1^{\hspace{-.5pt}\rule[2pt]{3pt}{.3pt}
\hspace{-1.5pt}\rule[.7pt]{.3pt}{3pt}}}
\begin{document}
\begin{abstract} We prove that the normalizer of any diffuse
amenable subalgebra of a free group factor $L(\Bbb F_r)$ generates
an amenable von Neumann subalgebra. Moreover, any II$_1$ factor of
the form $Q \vt L(\Bbb F_r) $, with $Q$ an arbitrary subfactor of a
tensor product of free group factors, has no Cartan subalgebras. We
also prove that if a free ergodic measure preserving action of a
free group $\Bbb F_r$, $2\leq r \leq \infty$, on a probability space
$(X,\mu)$ is profinite then the group measure space factor
$L^\infty(X)\rtimes \Bbb F_r$ has unique Cartan subalgebra, up to
unitary conjugacy.
\end{abstract}
\maketitle
\section{Introduction}
A celebrated theorem of Connes (\cite{connes:cls}) shows that all
amenable II$_1$ factors are isomorphic to the approximately finite
dimensional (AFD) II$_1$ factor $R$ of Murray and von Neumann
(\cite{mvn}). In particular, all II$_1$ group factors $L(\Gamma)$
associated with ICC (infinite conjugacy class) amenable groups
$\Gamma$, and all group measure space II$_1$ factors
$L^\infty(X)\rtimes \Gamma$ arising from free ergodic measure
preserving actions of countable amenable groups $\Gamma$ on a
probability space $\Gamma \curvearrowright X$, are isomorphic to
$R$. Moreover, by \cite{cfw}, any decomposition of $R$ as a group
measure space algebra is unique, i.e.\ if $R=L^\infty(X_i)\rtimes
\Gamma_i$, for some free ergodic measure preserving actions
$\Gamma_i \curvearrowright X_i$, $i=1,2$, then there exists an
automorphism of $R$ taking $L^\infty(X_1)$ onto $L^\infty(X_2)$. In
fact, any two Cartan subalgebras of $R$ are conjugate by an
automorphism of $R$.

Recall in this respect that a {\it Cartan subalgebra} $A$ in a
II$_1$ factor $M$ is a maximal abelian $^*$-subalgebra $A\subset M$
with normalizer $\mathcal N_M(A)=\{u\in \mathcal U(A) \mid
uAu^*=A\}$ generating $M$  (\cite{dixmier}, \cite{feldman-moore}).
Its presence amounts to realizing $M$ as a generalized (twisted)
version of the group measure space construction, for an action
$\Gamma \curvearrowright X$ and a 2-cocycle, with $A=L^\infty (X)$.
Decomposing factors this way is important, especially if one can
show uniqueness of their Cartan subalgebras, because then the
classification of the factors reduces to the classification of the
corresponding actions $\Gamma \curvearrowright X$ up to orbit
equivalence (\cite{feldman-moore}). But beyond the amenable case,
very little is known about uniqueness, or possible non-existence, of
Cartan subalgebras in group factors, or other factors that are
\textit{a priori} constructed in different ways than as group measure space
algebras.

We investigate in this paper Cartan decomposition properties for a
class of non-amenable II$_1$ factors that are in some sense
``closest to being amenable''. Thus, we consider factors $M$
satisfying the {\it complete metric approximation property} ({\it
c.m.a.p.}) of Haagerup (\cite{haagerup:map}), which requires
existence of normal, finite rank, completely bounded (cb) maps
$\phi_n\colon M\rightarrow M$, such that $\|\phi_n\|_{\cb}\leq 1$
and $\lim \|\phi_n(x)-x\|_2=0$, $\forall x\in M$, where
$\|\cdot\|_2$ denotes the Hilbert norm given by the trace of $M$
(note that if $\phi_n$ could be taken unital, $M$ would follow
amenable). This is same as saying that the Cowling-Haagerup constant
$\Lambda_{\cb}(M)$ equals $1$ (see \cite{cowling-haagerup}). The
prototype non-amenable c.m.a.p. factors are the free group factors
$L(\F_r)$, $2 \leq r \leq \infty$ (\cite{haagerup:map}). Like
amenability, the c.m.a.p. passes to subfactors and is well behaved
to inductive limits and tensor products.

We in fact restrict our attention to c.m.a.p. factors of the form
$M=Q \rtimes \F_r$, and to subfactors $N$ of such $M$. Our aim is to
locate all (or prove possible absence of) diffuse AFD subalgebras
$P\subset N$ whose normalizer $\Nor_N(P)$ generates $N$. Our
general result along these lines shows:

\begin{thm*} Let $\F_r \curvearrowright Q$ be an action of a free
group on a finite von Neumann algebra. Assume $M=Q\rtimes \F_r$ has
the complete metric approximation property. If $P\subset M$ is a
diffuse amenable subalgebra and $N$ denotes the von Neumann algebra
generated by its normalizer $\Nor_M(P)$, then either $N$ is
amenable relative to $Q$ inside $M$, or $P$ can be embedded into $Q$
inside $M$.
\end{thm*}

The amenability property of a von Neumann subalgebra $N\subset M$
relative to another von Neumann subalgebra $Q\subset M$ is rather
self-explanatory: it requires existence of a norm one projection
from the basic construction algebra of the inclusion $Q\subset M$
onto $N$ (see Definition \ref{defn:relamen}). The ``embeddability
of a subalgebra $P\subset M$ into another subalgebra $Q\subset M$
inside an ambient factor $M$'' is in the sense of \cite{popa:strong}
(see Definition \ref{defn:embed} below), and roughly means that
$P$ can be conjugated into $Q$ via a unitary element of $M$.

We mention three applications of the Theorem, each corresponding to
a particular choice of $\F_r \curvearrowright Q$ and solving well
known problems. Thus, taking $Q=\C$, we get:

\begin{cor*} The normalizer of any diffuse amenable
subalgebra $P$ of a free group factor $L(\F_r)$ generates an
amenable $($thus AFD by \cite{connes:cls}$)$ von Neumann algebra.
\end{cor*}

If we take $Q$ to be an arbitrary finite factor with
$\Lambda_{\cb}(Q)=1$ and let $\F_r$ act trivially on it, then
$M=Q\vt L(\F_r)$, $\Lambda_{\cb}(M)=1$ and the Theorem implies:

\begin{cor*} If $Q$ is a $\mathrm{II}_1$
factor with the complete metric approximation
property then $Q\vt L(\F_r)$ does not have Cartan subalgebras.
Moreover, if $N \subset Q\vt L(\F_r)$ is a subfactor of finite index
\cite{jones}, then $N$ does not have Cartan subalgebras either.
\end{cor*}

This shows in particular that any factor of the form $L(\F_r) \vt
R$, $L(\F_{r_1}) \vt L(\F_{r_2}) \vt\cdots$, and more generally any
subfactor of finite index of such a factor, has no Cartan
decomposition. Besides $Q=R, L(\F_r)$, other examples of factors
with $\Lambda_{\cb}(Q)=1$ are the group factors $L(\Gamma)$
corresponding to ICC discrete subgroups $\Gamma$ of
$\mathrm{SO}(1,n)$ and $\mathrm{SU}(1,n)$
(\cite{decanniere-haagerup,cowling-haagerup}), as well as any
subfactor of a tensor product of such factors. None of the factors
covered by Corollary 2 were known until now not to have Cartan
decomposition.

Finally, if we take $\F_r \curvearrowright X$ to be a {\it
profinite} measure preserving action on a probability measure space
$(X,\mu)$, i.e.\ an action with the property that $L^\infty(X)$ is a
limit of an increasing sequence of $\F_r$-invariant finite
dimensional subalgebras $Q_n$ of $L^\infty(X)$, then
$M=L^\infty(X)\rtimes \F_r$ is an increasing limit of the algebras
$Q_n \rtimes \F_r$, each one of which is an amplification of
$L(\F_r)$. Since c.m.a.p. behaves well to amplifications and
inductive limits, it follows that $M$ has c.m.a.p., so by applying
the Theorem and (A.1 in \cite{popa:betti}) we get:

\begin{cor*} If $\F_r\curvearrowright X$ is a free ergodic measure
preserving profinite action, then $L^\infty(X)$ is the unique Cartan
subalgebra of the $\mathrm{II}_1$-factor $L^\infty(X)\rtimes\F_r$,
up to unitary conjugacy.
\end{cor*}

The above corollary produces the first examples of non-amenable
II$_1$ factors with {\it all} Cartan subalgebras unitary conjugate.
Indeed, the ``unique Cartan decomposition'' results in
\cite{popa:betti,popa:strong,ipp} only showed conjugacy of Cartan
subalgebras satisfying certain properties. This was still enough for
differentiating factors of the form $L^\infty({\mathbb T}^2) \rtimes
\F_r$ and calculating their fundamental group in \cite{popa:betti},
by using \cite{gaboriau}. Similarly here, when combined with
Gaboriau's results, Corollary 3 shows that any factor
$L^\infty(X)\rtimes \F_r$, $2\leq r < \infty$, arising from a free
ergodic profinite action $\F_r \curvearrowright X$, has trivial
fundamental group. Also, if $\F_s \curvearrowright X$ is another
such action, with $r < s \leq \infty$, then $L^\infty(X)\rtimes \F_r
\not\simeq L^\infty(Y)\rtimes \F_s$. It can be shown that the
factors considered in \cite{popa:betti,popa:strong,ipp} cannot even
be embedded into the factors arising from profinite actions of free
groups. Note that the uniqueness of the Cartan subalgebras of the
AFD factor $R$ is up to conjugacy by automorphisms (\cite{cfw}),
but not up to unitary conjugacy, i.e.\ up to conjugacy by inner
automorphisms. Indeed, by \cite{feldman-moore} there exist
uncountably many non-unitary conjugate Cartan subalgebras in $R$.
Finally, note that Connes and Jones constructed examples of II$_1$
factors $M$ with two Cartan subalgebras that are not conjugate by
automorphisms of $M$ (\cite{connes-jones}).

Corollary 1 strengthens two well known in-decomposability properties
of free group factors: Voiculescu's result in \cite{voiculescu},
showing that $L(\F_r)$ has no Cartan subalgebras, which in fact
exhibited the first examples of factors with no Cartan
decomposition; and the first named author's result in
\cite{ozawa:solid}, showing that the commutant in $L(\F_r)$ of any
diffuse subalgebra must be amenable ($L(\F_r)$ are {\it solid}),
which itself strengthened the in-decomposability of $L(\F_r)$ into
tensor product of II$_1$ factors (\emph{primeness} of free group
factors) in \cite{ge}.

One should point out that Connes already constructed in
\cite{connes:chi} a factor $N$ that does not admit a ``classic''
group measure space decomposition $L^\infty(X)\rtimes \Gamma$. His
factor $N$ is defined as the fixed point algebra of an appropriate
finite group of automorphisms of $M=R \vt L(\F_r)$. But it was left
open whether $N$ cannot be obtained as a generalized group measure
space factor either, i.e.\ whether it does not have Cartan
decomposition. Corollary 2 shows that indeed it doesn't.

The proof of the Theorem follows a ``deformation/rigidity''
strategy, being inspired by arguments in \cite{popa:indec} and
\cite{popa:betti}. The proof is in two parts:

First we show that if a factor $M$ has the complete metric
approximation property then given any AFD  subalgebra $P\subset M$
the action (by conjugation) of the normalizer $\Nor_{M}(P)$ on
it satisfies a certain {\it weak compactness} property. This
essentially means $L^2(P)$ is a limit of finite dimensional
subspaces that are almost invariant to both the left multiplication
by elements in $P$ and to the $\Ad\Nor_{M}(P)$ action, in the
Hilbert-Schmidt norm (see Theorem \ref{nssg}). Note that this
implies wreath product factors $M=B^\Gamma \rtimes \Gamma$,
with $\Gamma$ non-amenable and $B\neq \C$, can never have the
complete metric approximation property. In particular
$\Lambda_{\cb}(\Lambda\wr\G)>1$, $\forall \Lambda \neq 1$, a fact
that was open until now.

For the second part, assume for simplicity $M=L(\F_r)$. Let
$P\subset M$ be diffuse AFD, $N=\Nor_M(P)''$. Taking $\eta\in
HS(L^2(M))\simeq L^2(M)\vt L^2(\bar{M})$ to be F{\o}lner-type
elements as given by the first part of the proof and $\alpha_t$ the
``malleable deformation'' of $L(\F_r)*L(\F_r)$ in
\cite{popa:bernoulli,popa:indec}, it follows that for $t$ small the
elements $(\alpha_t\otimes1)(\eta)\in L^2(M*M)\vt L^2(\bar{M})$
are still ``almost invariant,'' in the above sense.
This is used to prove that $L^2(N)$
is weakly contained in a multiple of the {\it coarse} bimodule
$L^2(M)\vt L^2(M)$, thus showing $N$ is AFD by the characterizations
of amenability in \cite{connes:cls}. This argument is the subject of
Theorem \ref{crossedprod} in the text.

We recall in Section \ref{sec:prelim} of the paper a number of known
results needed in the proofs, for the reader's convenience.
This includes a discussion of relative amenability (\ref{sec:relamen}),
intertwining lemmas for subalgebras (\ref{sec:bimodule}) and several
facts on the complete metric approximation property (\ref{sec:cmap}).
We mention that in the last Section of the paper we prove that for each
$2\leq r \leq \infty$ there exist uncountably many non orbit equivalent
profinite actions $\F_r\curvearrowright X$, which by Corollary 3 provide
uncountably many non-isomorphic factors $L^\infty(X)\rtimes \F_r$
as well (see \ref{cor:uncappfree}).
\section{Preliminaries}\label{sec:prelim}
\subsection{Finite von Neumann algebras}\label{sec:finitevna}
We fix conventions for (semi-)finite von Neumann algebras, but
before that we note that the symbol ``$\Lim$'' will be used for a
state on $\ell^\infty(\N)$, or more generally on $\ell^\infty(I)$
with $I$ directed, which extends the ordinary limit, and that the
abbreviation ``u.c.p.''\ stands for ``unital completely positive.''
We say a map is \emph{normal} if it is ultraweakly continuous.
Whenever a \emph{finite} von Neumann algebra $M$ is being
considered, it comes equipped with a distinguished faithful normal
tracial state, denoted by $\tau$. Any group action on a finite von
Neumann algebra is assumed to preserve the tracial state $\tau$. If
$M=L(\G)$ is a group von Neumann algebra, then the tracial state
$\tau$ is given by $\tau(x)=\ip{x\delta_1,\delta_1}$ for $x\in
L(\G)$. Any von Neumann subalgebra $P\subset M$ is assumed to contain
the unit of $M$ and inherits the tracial state $\tau$ from $M$.
The unique $\tau$-preserving conditional expectation from $M$ onto $P$
is denoted by $E_P$. We denote by $\Zt(M)$ the center of $M$;
by $\U(M)$ the group of unitary elements in $M$; and by
\[
\Nor_M(P)=\{ u\in\U(M) : (\Ad u)(P)=P\}
\]
the normalizing group of $P$ in $M$, where $(\Ad u)(x)=uxu^*$.
A maximal abelian von Neumann subalgebra $A\subset M$ satisfying
$\Nor_M(A)''=M$ is called a {\it Cartan subalgebra}.
We note that if $\G\curvearrowright X$ is an ergodic
essentially-free probability-measure-preserving action, then
$A=L^\infty(X)$ is a Cartan subalgebra in the crossed product
$L^\infty(X)\rtimes\G$. (See \cite{feldman-moore}.)

We refer the reader to the section IX.2 of \cite{takesaki:II} for
the details of the following facts on noncommutative $L^p$-spaces.
Let $\NN$ be a semi-finite von Neumann algebra with a faithful
normal semi-finite trace $\Tr$. For $1\le p<\infty$, we define the
$L^p$-norm on $\NN$ by $\|x\|_p=\Tr(|x|^p)^{1/p}$. By completing
$\{x\in\NN : \|x\|_p<\infty\}$ with respect to the $L^p$-norm, we
obtain a Banach space $L^p(\NN)$. We only need $L^1(\NN)$,
$L^2(\NN)$ and $L^\infty(\NN)=\NN$. The trace $\Tr$ extends to a
contractive linear functional on $L^1(\NN)$. We occasionally write
$\widehat{x}$ for $x\in\NN$ when viewed as an element in $L^2(\NN)$.
For any $1\le p,q,r\le\infty$ with $1/p+1/q=1/r$, there is a natural
product map
\[
L^p(\NN)\times L^q(\NN)\ni(x,y)\mapsto xy\in L^r(\NN)
\]
which satisfies $\|xy\|_r\le\|x\|_p\|y\|_q$
for any $x$ and $y$.
The Banach space $L^1(\NN)$ is identified with
the predual of $\NN$ under the duality
$L^1(\NN)\times\NN\ni(\zeta,x)\mapsto\Tr(\zeta x)\in\C$.
The Banach space $L^2(\NN)$ is identified with the GNS-Hilbert
space of $(\NN,\Tr)$.
Elements in $L^p(\NN)$ can be regarded as closed operators
on $L^2(\NN)$ which are affiliated with $\NN$ and hence
in addition to the above-mentioned product, there are
well-defined notion of positivity, square root, etc.
We will use many times the generalized
Powers--St{\o}rmer inequality (Theorem XI.1.2 in \cite{takesaki:II}):
\begin{equation}\label{eq:ps}
\|\eta-\zeta\|_2^2\le\|\eta^2-\zeta^2\|_1\le\|\eta+\zeta\|_2\|\eta-\zeta\|_2
\end{equation}
for every $\eta,\zeta\in L^2(\NN)_+$.
The Hilbert space $L^2(\NN)$ is an $\NN$-bimodule
such that $\ip{x\xi y,\eta}=\Tr(x\xi y\eta^*)$
for $\xi,\eta\in L^2(\NN)$ and $x,y\in\NN$.
We recall that this gives the canonical identification
between the commutant $\NN'$ of $\NN$ in $\B(L^2(\NN))$ and
the opposite von Neumann algebra
$\NN^{\op}=\{x^{\op} :x\in\NN\}$ of $\NN$.
Moreover, the opposite von Neumann algebra $\NN^{\op}$
is $*$-isomorphic to the complex conjugate von Neumann
algebra $\bar{\NN}=\{\bar{x} : x\in\NN\}$ of $\NN$ under
the $*$-isomorphism $x^{\op}\mapsto\bar{x}^*$.

Whenever $\NN_0\subset\NN$ is a von Neumann subalgebra
such that the restriction of $\Tr$ to $\NN_0$ is
still semi-finite, we identify $L^p(\NN_0)$ with
the corresponding subspace of $L^p(\NN)$.
Anticipating a later use, we consider the tensor product
von Neumann algebra $(\NN\vt M,\Tr\otimes\tau)$
of a semi-finite von Neumann algebra $(\NN,\Tr)$ and
a finite von Neumann algebra $(M,\tau)$.
Then, $\NN\cong\NN\vt\C1\subset\NN\vt M$ and
the restriction of $\Tr\otimes\tau$ to $\NN$ is $\Tr$.
Moreover, the conditional expectation
$\id\otimes\tau\colon\NN\vt M\to\NN$ extends
to a contraction from $L^1(\NN\vt M)\to L^1(\NN)$.

Let $Q\subset M$ be finite von Neumann algebras.
Then, the conditional expectation $E_Q$ can be viewed
as the orthogonal projection $e_Q$ from $L^2(M)$ onto
$L^2(Q)\subset L^2(M)$.
It satisfies $e_Qxe_Q=E_Q(x)e_Q$ for every $x\in M$.
The \emph{basic construction} $\ip{M,e_Q}$ is the
von Neumann subalgebra of $\B(L^2(M))$ generated by
$M$ and $e_Q$. We note that $\ip{M,e_Q}$ coincides
with the commutant of the right $Q$-action in $\B(L^2(M))$.
The linear span of $\{ xe_Qy : x,y\in M\}$ is an ultraweakly
dense $*$-subalgebra in $\ip{M,e_Q}$ and the basic construction
$\ip{M,e_Q}$ comes together with the faithful normal semi-finite
trace $\Tr$ such that $\Tr(xe_Qy)=\tau(xy)$.
See Section 1.3 in \cite{popa:betti} for more information on
the basic construction.
\subsection{Relative amenability}\label{sec:relamen}
We adapt here Connes's characterization of amenable von Neumann
algebras to the relative situation. Recall that for von Neumann
algebras $N\subset\NN$, a state $\p$ on $\NN$ is said to be
\emph{$N$-central} if $\p\circ\Ad(u)=\p$ for any $u\in\U(N)$, or
equivalently if $\p(ax)=\p(xa)$ for all $a\in N$ and $x\in\NN$.

\begin{thm}\label{relativeamenability}
Let $Q, N\subset M$ be finite von Neumann algebras.
Then, the following are equivalent.
\begin{enumerate}
\item\label{con:relamen}
There exists a $N$-central state $\p$
on $\ip{M,e_Q}$ such that $\p|_M=\tau$.
\item\label{con:relamen2}
There exists a $N$-central state $\p$
on $\ip{M,e_Q}$ such that $\p$ is normal on $M$
and faithful on $\Zt(N'\cap M)$.
\item\label{con:condexp}
There exists a conditional expectation $\Phi$ from $\ip{M,e_Q}$ onto $N$
such that $\Phi|_M=E_N$.
\item\label{con:ell2global}
There exists a net $(\xi_n)$ in $L^2\ip{M,e_Q}$ such
that $\lim_n\ip{x\xi_n,\xi_n}=\tau(x)$ for every $x\in M$
and that $\lim\|[u,\xi_n]\|_2=0$ for every $u\in N$.
\end{enumerate}
\end{thm}
\begin{defn}\label{defn:relamen}
Let $Q,N\subset M$ be finite von Neumann algebras.
We say $N$ is \emph{amenable relative to $Q$ inside $M$},
denoted by $N\lessdot_MQ$, if any of the conditions in
Theorem \ref{relativeamenability} holds.
We say $Q$ is \emph{co-amenable in $M$} if $M\lessdot_MQ$
(cf.\ \cite{popa:corresp,delaroche}).
\end{defn}
\begin{proof}[Proof of Theorem \ref{relativeamenability}]
The proof follows a standard recipe
of the theory (cf.\ \cite{connes:cls,haagerup:dec,popa:corresp}).
The implication $(\ref{con:relamen})\Rightarrow(\ref{con:relamen2})$ is obvious.
To prove the converse, assume the condition $(\ref{con:relamen2})$.
Then, there exists $b\in L^1(M)_+$ such that $\p(x)=\tau(bx)$ for $x\in M$.
Since $\p$ is $N$-central, one has $ubu^*=b$ for all $u\in\U(N)$,
i.e.\ $b\in L^1(N'\cap M)$.
We consider the directed set $I$ of finite subsets of $\U(N'\cap M)$.
For each element $i=\{u_1,\ldots,u_n\}\in I$ and $m\in\N$, we define
$b_i=n^{-1}\sum u_kbu_k^*\in L^1(N'\cap M)_+$,
$c_{i,m}=\chi_{(1/m,\infty)}(b_i)b_i^{-1/2}\in N'\cap M$
and
\[
\psi_{i,m}(x)=\frac{1}{n}\sum_{k=1}^n\p(u_k^*c_{i,m}xc_{i,m}u_k)
\]
for $x\in\ip{M,e_Q}$.
Since $c_{i,m}u_k\in N'\cap M$, the positive linear functionals
$\psi_{i,m}$ are still $N$-central and
$\psi_{i,m}(x)=\tau(\chi_{(1/m,\infty)}(b_i)x)$ for $x\in M$.
We note that
\[
\lim_i\lim_m \chi_{(1/m,\infty)}(b_i)
=\lim_i s(b_i)
=\lim_i\bigvee s(u_kbu_k^*)=z,
\]
where $s(\,\cdot\,)$ means the support projection and
$z$ is the central support projection of $b$ in $N'\cap M$.
Since $\p(z^\perp)=\tau(bz^\perp)=0$ and $\p$ is faithful on
$\Zt(N'\cap M)$, one has $z=1$.
Hence, the state $\psi=\Lim_i\Lim_m\psi_{i,m}$ on $\ip{M,e_Q}$
is $N$-central and satisfies $\psi|_M=\tau$. This proves $(\ref{con:relamen})$.

We prove $(\ref{con:relamen})\Rightarrow(\ref{con:ell2global})$:
Let a $N$-central state $\p$ on $\ip{M,e_Q}$ be given
such that $\p|_M=\tau$.
Take a net $(\zeta_n)$ of positive norm-one elements in $L^1\ip{M,e_Q}$
such that $\Tr(\zeta_n\,\cdot\,)$ converges to $\p$ pointwise.
Then, for every $x\in\ip{M,e_Q}$ and $u\in\U(N)$, one has
\[
\lim_n\Tr((\zeta_n-\Ad(u)\zeta_n)x)=\p(x)-\p(\Ad(u^*)(x))=0
\]
by assumption.
It follows that for every $u\in\U(N)$,
the net $\zeta_n-\Ad(u)(\zeta_n)$ in $L^1\ip{M,e_Q}$
converges to zero in the weak-topology.
By the Hahn-Banach separation theorem,
one may assume, by passing to convex combinations,
that it converges to zero in norm.
Thus, $\|[u,\zeta_n]\|_1\to0$ for every $u\in\U(N)$.
By (\ref{eq:ps}), if we define
$\xi_n=\zeta_n^{1/2}\in L^2\ip{M,e_Q}$,
then one has
$\|[u,\xi_n]\|_2\to0$ for every $u\in\U(N)$.
Moreover, for any $x\in M$,
\[
\lim_n\ip{x\xi_n,\xi_n}=\lim_n\Tr(\zeta_n x)=\p(x)=\tau(x).
\]

We prove $(\ref{con:ell2global})\Rightarrow(\ref{con:condexp})$:
Since
\[
|\p(bcyz)|=|\p(cyzb)|
\le\p(cyy^*c^*)^{1/2}\p(b^*z^*zb)^{1/2}
\le\|b\|_2\|c\|_2\|y\|\|z\|
\]
for every $b,c\in N$ and $y,z\in\ip{M,e_Q}$,
one has $|\p(ax)|\le\|a\|_1\|x\|$ for
every $a\in N$ and $x\in\ip{M,e_Q}$.
Hence, for every $x\in\ip{M,e_Q}$, we may
define $\Phi(x)\in N=L^1(N)^*$ by the duality
$\tau(a\Phi(x))=\p(ax)$ for all $a\in N$.
It is clear that $\Phi$ is a conditional expectation onto $N$
such that $\Phi|_M=E_N$.

We prove $(\ref{con:condexp})\Rightarrow(\ref{con:relamen})$:
If there is a conditional expectation $\Phi$ from $\ip{M,e_Q}$ onto $N$
such that $\Phi|_M=E_N$, then $\p=\tau\circ\Phi$ is
an $N$-central state such that $\p|_M=\tau$.
\end{proof}
Let $N_0\subset M$ be a von Neumann subalgebra whose unit $e$
does not coincide with the unit of $M$.
We say $N_0$ is amenable relative to $Q$ inside $M$,
denoted by $N_0\lessdot_MQ$, if $N_0+\C(1-e)\lessdot_MQ$.
We observe that $N_0\lessdot_MQ$ if and only if
there exists an $N_0$-central state $\p$ on $e\ip{M,e_Q}e$
such that $\p(exe)=\tau(exe)/\tau(e)$ for $x\in M$.
\begin{cor}\label{cor:relative}
Let $Q_1,\ldots,Q_k,N\subset M$ be finite von Neumann algebras
and $\GG\subset\U(N)$ be a subgroup such that $\GG''=N$.
Assume that for every non-zero projection $p\in\Zt(N'\cap M)$,
there exists a net $(\xi_n)$ of vectors in
a multiple of $\bigoplus_{j=1}^k L^2\ip{M,e_{Q_j}}$ such that
\begin{enumerate}
\item $\limsup\|x\xi_n\|_2\le\|x\|_2$ for all $x\in M$;
\item $\liminf\|p\xi_n\|_2>0$; and
\item $\lim\|[u,\xi_n]\|_2=0$ for every $u\in\GG$.
\end{enumerate}
Then, there exist projections $p_1,\ldots,p_k\in\Zt(N'\cap M)$
such that $\sum_{j=1}^k p_j=1$ and $Np_j\lessdot_MQ_j$ for every $j$.
\end{cor}
\begin{proof}
We observe that if there exists an increasing net $(e_i)_i$ of
projections in $\Zt(N'\cap M)$ such that $Ne_i\lessdot_MQ$ for all $i$,
then $Ne\lessdot_MQ$ for $e=\sup e_i$.
Hence, by Zorn's lemma, there is a maximal $k$-tuple $(p_1,\ldots,p_k)$
of projections in $\Zt(N'\cap M)$ such that $\sum_j p_j\le1$
and $Np_j\lessdot_MQ_j$ for every $j$. We prove that $\sum_j p_j=1$.
Suppose by contradiction that $p=1-\sum_j p_j\neq0$,
and take a net $(\xi_n)$ as in the statement of the corollary.
We may assume that all $\xi_n$'s are in a multiple of $L^2\ip{M,e_{Q_j}}$
for some fixed $j\in\{1,\ldots,k\}$.
We define a state $\psi$ on $\ip{M,e_{Q_j}}$ by
\[
\psi(x)=\Lim_n \|p\xi_n\|_2^{-2}\ip{xp\xi_n,p\xi_n}
\]
for $x\in\ip{M,e_{Q_j}}$.
It is not hard to see that $\psi(p)=1$,
$\psi\circ\Ad(u)=\psi$ for every $u\in\GG$ and
$\psi(x^*x)\le (\liminf\|p\xi_n\|)^{-2}\|xp\|_2^2$ for every $x\in M$.
It follows that $\psi|_M$ is normal and $\psi$ is $N$-central.
Let $q$ be the minimal projection in $\Zt(N'\cap M)$ such that $\psi(q)=1$.
We finish the proof by showing $Nr\lessdot_MQ_j$ for $r=p_j+q$
(which gives the desired contradiction to maximality).
Since $Np_j\lessdot_MQ_j$, there is an $Np_j$-central state $\p$
on $p_j\ip{M,e_{Q_j}}p_j$ such that $\p(p_jxp_j)=\tau(p_jxp_j)/\tau(p_j)$ for $x\in M$.
We fix a state extension $\tilde{\tau}$ of $\tau$ on $\ip{M,e_{Q_j}}$ and
define a state $\tilde{\p}$ on $\ip{M,e_{Q_j}}$ by
\[
\tilde{\p}(x)=\tau(p_j)\p(p_jxp_j)+\tau(q)\psi(qxq)+\tilde{\tau}((1-r)x(1-r))
\]
for $x\in\ip{M,e_{Q_j}}$.
The state $\tilde{\p}$ is $(Nr+\C(1-r))$-central,
normal on $M$ and faithful on $\Zt((Nr+\C(1-r))'\cap M)=\Zt(N'\cap M)r+\Zt(M)(1-r)$.
Hence Theorem~\ref{relativeamenability} implies $Nr\lessdot_MQ_j$.
\end{proof}
Compare the following result
with \cite{popa:corresp} and \cite{delaroche}.

\begin{prop}\label{proprelamen}
Let $P,Q,N\subset M$ be finite von Neumann algebras.
Then, the following are true.
\begin{enumerate}
\item
Suppose that $M=Q\rtimes\G$ is the crossed product of
$Q$ by a group $\G$.
Then, $L(\G)\lessdot_MQ$ if and only if $\G$ is amenable.
\item
Suppose that $Q$ is AFD. Then, $P\lessdot_MQ$
if and only if $P$ is AFD.
\item
If $N\lessdot_MP$ and $P\lessdot_MQ$, then $N\lessdot_MQ$.
\end{enumerate}
\end{prop}

\begin{proof}
Denote by $\lambda_g$ the unitary element in $M$ which implements
the action of $g\in\G$. Since $e_Q\lambda(g)e_Q=0$ for
$g\in\G\setminus\{1\}$, the projections
$\{\lambda_ge_Q\lambda_g^*:g\in\G\}$ are mutually orthogonal and
generate an isomorphic copy of $\ell^\infty(\G)$ in $\ip{M,e_Q}$.
Hence, if there exists an $L(\G)$-central state on $\ip{M,e_Q}$, then its
restriction to $\ell^\infty(\G)$ becomes a $\G$-invariant mean. This
proves the ``only if'' part of the assertion (1). The ``if'' part is
trivial. The assertion (2) easily follows from the fact that
$\ip{M,e_Q}$ is injective if (and only if) $Q$ is AFD
(\cite{connes:cls}).

Let us finally prove (3). Fix a conditional expectation $\Phi$ from
$\ip{M,e_Q}$ onto $P$ such that $\Phi|_M=E_P$. For $\xi=\sum_{i=1}^m
a_i\otimes b_i\in M\otimes M$, we denote
\[
\|\xi\|_2 =\|\sum_{i=1}^m a_i e_P b_i\|_{L^2\ip{M,e_P}}
 =\bigl(\sum_{i,j}\tau(b_i^*E_P(a_i^*a_j)b_j)\bigr)^{1/2}.
\]
For $\xi=\sum_{i=1}^m a_i\otimes b_i$ and $\eta=\sum_{j=1}^n c_j\otimes d_j$
in $M\otimes M$,
we define a linear functional $\p_{\eta,\xi}$ on $\ip{M,e_Q}$ by
\[
\p_{\eta,\xi}(x)
=\sum_{i,j}\tau(b_i^*\Phi(a_i^*xc_j)d_j).
\]
We claim that $\|\p_{\eta,\xi}\|\le\|\eta\|_2\|\xi\|_2$.
Indeed, if $\Phi(x)=V^*\pi(x)V$ is a Stinespring dilation,
then one has
\[
\p_{\eta,\xi}(x)=\ip{\pi(x)\sum_j\pi(c_j)Vd_j\widehat{1}_P,
\sum_i\pi(a_i)Vb_i\widehat{1}_P}
\]
and $\|\sum_i\pi(a_i)Vb_i\widehat{1}_P\|=\|\xi\|_2$
and likewise for $\eta$.
It follows that $\p_{\eta,\xi}$ is defined for
$\xi,\eta\in L^2\ip{M,e_P}$ in such a way that
$\|\p_{\eta,\xi}\|\le\|\eta\|_2\|\xi\|_2$.
Now take a net of unit vectors $(\xi_n)$ in $L^2\ip{M,e_P}$
satisfying the condition \ref{con:ell2global}
in Theorem \ref{relativeamenability}, and let
$\p=\Lim\p_{\xi_n,\xi_n}$ be the state on $\ip{M,e_Q}$.
Then, one has
\[
\p\circ\Ad(u)=\Lim_n\p_{\Ad(u)(\xi_n),\Ad(u)(\xi_n)}
=\Lim_n\p_{\xi_n,\xi_n}=\p
\]
for all $u\in\U(N)$ and
\[
\p(x)=\Lim_n\ip{x\xi_n,\xi_n}_{L^2\ip{N,e_P}}=\tau(x)
\]
for all $x\in M$. This proves that $N\lessdot_MQ$.
\end{proof}
\subsection{Intertwining subalgebras inside II$_1$ factors}\label{sec:bimodule}
We extract from \cite{popa:betti,popa:strong} some results
which are needed later. The following are Theorem A.1 in
\cite{popa:betti} and its corollary (also, a particular case of 2.1
in \cite{popa:strong}).

\begin{thm}\label{finitebimodule}
Let $N$ be a finite von Neumann algebra and $P,Q\subset N$
be von Neumann subalgebras. Then, the following are equivalent.
\begin{enumerate}
\item
There exists a non-zero projection $e\in\ip{N,e_Q}$ with $\Tr(e)<\infty$
such that the ultraweakly closed convex hull of $\{ w^*ew : w\in\U(P)\}$
does not contain $0$.
\item
There exist non-zero projections
$p\in P$ and $q\in Q$, a normal $*$-homomorphism
$\theta\colon pPp\to qQq$ and a non-zero partial isometry $v\in N$
such that
\[
\forall x\in pPp\quad xv=v\theta(x)
\]
and $v^*v\in\theta(pPp)'\cap qNq$, $vv^*\in p(P'\cap N)p$.
\end{enumerate}
\end{thm}

\begin{defn}\label{defn:embed}
Let $P,Q\subset N$ be finite von Neumann algebras. Following
\cite{popa:strong}, we say that \emph{$P$ embeds into $Q$ inside
$N$}, and write $P\preceq_NQ$, if any of the conditions in
Theorem~\ref{finitebimodule} holds.
\end{defn}

Let $\phi$ be a $\tau$-preserving u.c.p.\ map on $N$.
Then, $\phi$ extends to a contraction $T_{\phi}$
on $L^2(N)$ by $T_{\phi}(\widehat{x})=\widehat{\phi(x)}$.
Suppose that $\phi|_Q=\id_Q$. Then, $\phi$
automatically satisfies $\phi(axb)=a\phi(x)b$
for any $a,b\in Q$ and $x\in N$.
It follows that $T_{\phi}\in\B(L^2(N))$ commutes with
the right action of $Q$, i.e.,
$T_{\phi}\in\ip{N,e_Q}$.
We say $\phi$ is \emph{compact over} $Q$ if $T_{\phi}$ belongs
to the ``compact ideal'' of $\ip{N,e_Q}$
(see Section 1.3.2 in \cite{popa:betti}).
If $\phi$ is compact over $Q$, then for any $\e>0$,
the spectral projection $e=\chi_{[\e,1]}(T_{\phi}^*T_{\phi})\in\ip{N,e_Q}$
has finite $\Tr(e)$ and
\[
\ip{w^*ew\widehat{1},\widehat{1}}_{L^2(N)}
\geq\ip{T_{\phi}^*T_{\phi}\widehat{w},\widehat{w}}_{L^2(N)}-\e
=\|\phi(w)\|_2^2-\e
\]
for all $w\in\U(P)$.
These observations imply the following corollary \cite{popa:betti}.

\begin{cor}\label{cor:cptq}
Let $P,Q\subset N$ be finite von Neumann algebras.
Suppose that $\phi$ is a $\tau$-preserving u.c.p.\ map
on $N$ such that $\phi|_Q=\id_Q$ and $\phi$ is compact over $Q$.
If $\inf\{\|\phi(w)\|_2 : w\in\U(P)\}>0$, then $P\preceq_NQ$.
\end{cor}

Finally, recall that A.1 in \cite{popa:betti} shows the following:
\begin{lem}\label{conjugatecartan}
Let $A$ and $B$ be maximal abelian $^*$-subalgebras of a type
$\mathrm{II}_1$-factor $N$ such that $\mathcal N_N(A)'', \mathcal
N_N(B)''$ are factors (i.e. $A, B$ are semiregular \cite{dixmier}).
If $A\preceq_N B$, then there exists $u\in\U(N)$ such that
$uAu^*=B$.
\end{lem}
\subsection{The complete metric approximation property}\label{sec:cmap}
Let $\G$ be a discrete group.
For a function $f$ on $\G$, we write $m_f$ for the multiplier
on $\C\G\subset L(\G)$ defined by $m_f(g)=fg$ for $g\in\C\G$.
We simply write $\|f\|_{\cb}$ for $\|m_f\|_{\cb}$ and call it
the Herz-Schur norm.
If $\|f\|_{\cb}$ is finite and $f(1)=1$,
then $m_f$ extends to a $\tau$-preserving
normal unital map on $L(\G)$.
We refer the reader to sections 5 and 6 in \cite{pisier}
for an account of Herz-Schur multipliers.

\begin{defn}\label{defn:cmap}
A discrete group $\G$ is \emph{weakly amenable} if there exist a
constant $C\geq1$ and a net $(f_n)$ of finitely supported functions
on $\G$ such that $\limsup\|f_n\|_{\cb}\le C$ and $f_n\to1$
pointwise. The Cowling-Haagerup constant $\Lambda_{\cb}(\G)$ of $\G$
is defined as the infimum of the constant $C$ for which a net
$(f_n)$ as above exists.

We say a von Neumann algebra $M$ has the \emph{$($weak$^*)$
completely bounded approximation property} if there exist a constant
$C\geq1$ and a net $(\phi_n)$ of normal
finite-rank maps on $M$ such that $\limsup\|\phi_n\|_{\cb}\le
C$ and $\|x-\phi_n(x)\|_2\to0$ for every $x\in M$. The
Cowling-Haagerup constant $\Lambda_{\cb}(M)$ of $M$ is defined as
the infimum of the constant $C$ for which a net $(\phi_n)$ as above
exists. Also, we say that $M$ has the \emph{$($weak$^*)$ complete
metric approximation property} ({\it c.m.a.p.}) if
$\Lambda_{\cb}(M)=1.$
\end{defn}

By routine perturbation arguments,
one may arrange $\phi_n$'s in the above definition to be unital
and trace-preserving when $M$ is finite.
We are interested here in
the case $\Lambda_{\cb}(M)=1$, i.e.\ when $M$ has the complete metric
approximation property. We summarize below some known results in
this direction. For part $(7)$, recall that an action
of a group $\Gamma$ on a finite von Neumann algebra
$P$ is {\it profinite} if there exists an increasing sequence
of $\Gamma$-invariant
finite dimensional von Neumann subalgebras $P_n \subset P$ that generate $P$.
Note that this implies $P$ is AFD. If $P=L^\infty(X)$ is abelian and $\Gamma
\curvearrowright P$ comes from a measure preserving action
$\Gamma \curvearrowright X$, then the profiniteness
of $\Gamma \curvearrowright P$ amounts to
existence of a sequence of
$\Gamma$-invariant finite partitions
of $X$ that generate the $\sigma$-algebra of measurable subsets of $X$.

\begin{thm}\label{Lambdacb}
\begin{enumerate}
\item
$\Lambda_{\cb}(L(\G))=\Lambda_{\cb}(\G)$ for any $\G$.
\item
If $\G$ is a discrete subgroup of $\mathrm{SO}(1,n)$ or of $\mathrm{SU}(1,n)$,
then $\Lambda_{\cb}(\G)=1$.
\item
If $\G$ acts properly on a finite-dimensional CAT(0) cubical complex, then $\Lambda_{\cb}(\G)=1$.
\item
If $\Lambda_{\cb}(\G_i)=1$ for $i=1,2$,
then $\Lambda_{\cb}(\G_1\times\G_2)=1$ and $\Lambda_{\cb}(\G_1\ast\G_2)=1$.
\item
If $N\subset M$ are finite von Neumann algebras, then
$\Lambda_{\cb}(N)\le\Lambda_{\cb}(M)$. Moreover, if $N, M$ are
factors and $[M:N]<\infty$, then $\Lambda_{\cb}(M)=\Lambda_{\cb}(N)$ and
$\Lambda_{\cb}(M^t)=\Lambda_{\cb}(M)$, $\forall t>0$.
\item
Let $M$ be a finite von Neumann algebra and $(M_n)$ be an increasing net of
von Neumann subalgebras of $M$ such that $M=(\bigcup M_n)''$.
Then, $\Lambda_{\cb}(M)=\sup\Lambda_{\cb}(M_n)$.
\item
If $P$ is a finite von Neumann algebra and $\G\curvearrowright P$ is
a profinite action,
then $\Lambda_{\cb}(P\rtimes\G)=\Lambda_{\cb}(\G)$.
\end{enumerate}
\end{thm}

The assertions (1), (2), (3) and (4) are respectively due to
\cite{cowling-haagerup},
\cite{decanniere-haagerup,cowling-haagerup}, \cite{guentner-higson}
and \cite{ricard-xu}. The rest are trivial. We will see in Corollary
\ref{criterion} that property (7) generalizes to compact actions of
groups $\Gamma$, and even to actions of $\Gamma$ that are ``weakly
compact'', in the sense of Definition \ref{defn:profinite}.

We prove in this paper a general property about normal amenable
subgroups of groups with $\Lambda_{\cb}$-constant equal to 1. While
this property is a consequence of Theorem~\ref{nssg} (via
$(\ref{prof3}) \Leftrightarrow(\ref{prof4})$ in
Proposition~\ref{profinite}), we give here a direct proof in
group-theoretic framework. To this end, note that if
$\Lambda\triangleleft\G$ is a normal subgroup then the semi-direct
product group $\Lambda\rtimes\G$ acts on $\Lambda$ by
$(a,g)b=agbg^{-1}$, for $(a,g)\in\Lambda\rtimes\G$ and
$b\in\Lambda$.
\begin{prop}
Suppose that $\G$ has an infinite normal amenable
subgroup $\Lambda\triangleleft\G$
and that $\Lambda_{\cb}(\G)=1$.
Then there exists a $\Lambda\rtimes\G$-invariant
mean on $\ell^\infty(\Lambda)$
(i.e., $\G$ is co-amenable in $\Lambda\rtimes\G$).
In particular, $\G$ is inner-amenable.
(See Section~\ref{sec:uncountable} for the definition of inner-amenability.)
\end{prop}
\begin{proof}
Let $f_n$ be a net of finitely supported functions such that
$\sup\|f_n\|_{\cb}=1$ and $f_n\to 1$ pointwise.
By the Bo\.zejko-Fendler theorem (Theorem 6.4 in \cite{pisier}),
there are Hilbert space vectors $\xi_n(a)$ and $\eta_n(b)$ of
norm at most one such that
$f_n(ab^{-1})=\ip{\eta_n(b),\xi_n(a)}$ for all $a,b\in\G$.
Then, for every $g\in\G$, one has
\begin{align*}
\lim_n\sup_{a\in\G}\|\xi_n(ga)-\xi_n(a)\|^2
&\le \lim_n\sup_{a\in\G}2\bigl(\|\xi_n(ga)-\eta_n(a)\|^2
 +\|\eta_n(a)-\xi_n(a)\|^2\bigr)\\
 &\le \lim_n 2(2-2\Re f_n(g)+2-2\Re f_n(1))=0,
\end{align*}
and similarly $\lim_n\sup_{b\in\G}\|\eta_n(gb)-\eta_n(b)\|=0$
for every $g\in\G$.
It follows that
\[
\lim_n\|f_n-f_n^g\|_{\cb}=0
\]
for every $g\in\G$, where $f_n^g\in\C\G$
is defined by $f_n^g(a)=f_n(gag^{-1})$.
Now since $\Lambda\triangleleft\G$ is amenable,
the trivial representation $\tau_0\colon C^*_{\mathrm{red}}(\Lambda)\to\C$
is continuous. We define a linear functional $\omega_n$
on $C^*_{\mathrm{red}}(\Lambda)$ by
$\omega_n=\tau_0\circ m_{f_n}|_{C^*_{\mathrm{red}}(\Lambda)}$.
Since $f_n$ is finitely supported, $\omega_n$ is ultraweakly
continuous on $L(\Lambda)$.
We note that
$\lim\omega_n(\lambda(a))=1$ for all $a\in\Lambda$ and
\[
\lim_n\|\omega_n-\omega_n\circ\Ad(g)\|\le\lim_n\|f_n-f_n^g\|_{\cb}=0
\]
for all $g\in\G$. Since $\|\omega_n\|\le1$ and $\lim\omega_n(1)=1$,
we have $\lim\|\omega_n-|\omega_n|\|=0$.
We view $|\omega_n|$ as an element in $L^1(L(\Lambda))$
(which is $L^1(\widehat{\Lambda})$ if $\Lambda$ is abelian)
and consider $\zeta_n=|\omega_n|^{1/2}\in L^2(L(\Lambda))=\ell^2(\Lambda)$.
Then, the net $(\zeta_n)$
satisfies $\lim_n\ip{\lambda(a)\zeta_n,\zeta_n}=1$
for all $a\in\Lambda$ and
$\lim_n\|\zeta_n-\Ad(g)(\zeta_n)\|_2=0$ for all $g\in\G$ by (\ref{eq:ps}).
Therefore, the state $\omega$ on
$\ell^\infty(\Lambda)\subset\B(\ell^2(\Lambda))$ defined by
\[
\omega(x)=\Lim_n\ip{x\zeta_n,\zeta_n}
=\Lim_n\sum_{a\in\Lambda}x(a)\zeta_n(a)^2
\]
is $\Lambda\rtimes\G$-invariant.
Since $\Lambda$ is infinite, the $\Lambda$-invariant mean $\omega$
is singular, i.e, $\zeta_n\to0$ weakly.
This implies inner-amenability of $\G$.
\end{proof}

Recall that the wreath product $\Lambda\wr\G$
of $\Lambda$ by $\G$ is defined as the semi-direct product
$(\bigoplus_\G \Lambda)\rtimes\G$ of $\bigoplus_\G \Lambda$ by
the shift action $\G\curvearrowright\bigoplus_\G \Lambda$.

\begin{cor}
If $\Lambda\neq\{1\}$ and $\G$ is non-amenable then
$\Lambda_{\cb}(\Lambda\wr\G)>1$, i.e.\ $L(\Lambda\wr\G)$ does not
have the complete metric approximation property.
\end{cor}
\begin{proof}
Suppose that $\Lambda_{\cb}(\Lambda\wr\G)=1$.
Passing to a subgroup if necessary, we may assume that
$\Lambda$ is cyclic (amenable).
We observe that the stabilizing subgroup
$\{ g\in\G : \sigma_g(a)=a\}$ of
any non-neutral element $a\in\bigoplus_\G\Lambda$
is finite.
It follows that there is a $\G$-equivariant u.c.p.\
map from $\ell^\infty(\G)$ into
$\ell^\infty((\bigoplus_\G\Lambda)\setminus\{1\})$.
It follows that there is no $\G$-invariant mean on
$(\bigoplus_\G\Lambda)\setminus\{1\}$ since $\G$ is non-amenable.
Hence, any $\G$-invariant mean on $\bigoplus_\G\Lambda$
has to be concentrated on $\{1\}$.
Such mean cannot be $(\bigoplus_\G\Lambda)$-invariant.
\end{proof}
\begin{rem}
Let $\G=(\Z/2\Z)\wr\F_2$.
Since $\Lambda_{\cb}$ is multiplicative (\cite{cowling-haagerup})
and satisfies $\Lambda_{\cb}(\Gamma) > 1$ (by 2.11 above),
the direct product $\bigoplus\G$ of infinitely many copies of $\G$ is not
weakly amenable, i.e. $\Lambda_{\cb}(\bigoplus\G)=\infty$.
It is plausible that $\G$ itself is not
weakly amenable. De Cornulier--Stalder--Valette (\cite{csv})
recently proved the surprising result that, despite
satisfying $\Lambda_{\cb}(\Gamma)>1$,
the group $\G$ (and hence $\bigoplus\G$) has
Haagerup's compact approximation property \cite{haagerup:map}.
Taken together, these results falsify one implication of the
so-called Cowling's conjecture,
which asserts that
Haagerup's compact approximation property for a group $\Gamma$
is equivalent to the condition $\Lambda_{\cb}(\Gamma)=1$.
But there are still no known examples of groups $\Gamma$ which satisfy
$\Lambda_{\cb}(\Gamma)=1$ but fail Haagerup's
compact approximation property.
\end{rem}
\section{Weakly compact actions}\label{sec:profinite}
We introduce in this section a new property for group actions,
weaker than compactness (thus weaker than profiniteness as well) and
closely related to the complete metric approximation property of the
corresponding crossed product algebras. The main result of
this section will show that if a II$_1$ factor $M$ has the c.m.a.p.
then given any maximal abelian subalgebra $A$ of $M$ the action on
$A$ of its normalizer, $\mathcal N_M(A)$, is weakly compact. Also,
if a group $\Gamma$ satisfies $\Lambda(\Gamma)=1$ and $\Gamma
\curvearrowright X$ is weakly compact, then $M=L^\infty(X)\rtimes
\Gamma$ has c.m.a.p.

\begin{defn}\label{defn:profinite}
Let $\sigma$ be an action of a group $\G$ on a finite von Neumann
algebra $P$. Recall that
$\sigma$ is called \emph{compact} if
$\sigma(\G)\subset\Aut(P)$ is pre-compact in the point-ultraweak topology.
We call the action $\sigma$ \emph{weakly compact} if
there exists a net $(\eta_n)$ of unit vectors in
$L^2(P\vt\bar{P})_+$ such that
\begin{enumerate}
\item $\|\eta_n-(v\otimes\bar{v})\eta_n\|_2\to0$ for every $v\in\U(P)$.
\item $\|\eta_n-(\sigma_g\otimes\bar{\sigma}_g)(\eta_n)\|_2\to0$ for every $g\in\G$.
\item $\ip{(x\otimes1)\eta_n,\eta_n}=\tau(x)=\ip{\eta_n,(1\otimes\bar{x})\eta_n}$
for every $x\in P$ and every $n$.
\end{enumerate}
Here, we consider the action $\sigma$ on $P$ as the corresponding
unitary representation on $L^2(P)$.
By the proof of Proposition \ref{profinite}, the condition (3) can be replaced
with a formally weaker condition
\begin{enumerate}
\item[$(3')$] $\ip{(x\otimes1)\eta_n,\eta_n}\to\tau(x)$ for every $x\in P$.
\end{enumerate}
\end{defn}

Weak compactness is manifestly weaker than profiniteness,
which is why in an initial version of this paper we called
it {\it weak profiniteness}. We are
very grateful to Adrian Ioana, who pointed out to
us that the condition is even weaker than compactness (cf.
$(\ref{prof2})\Rightarrow(\ref{prof3})$ below) and suggested
a change in terminology.

\begin{prop}\label{profinite}
Let $\sigma$ be an action of a group $\G$ on
a finite von Neumann algebra $P$ and consider
the following conditions.
\begin{enumerate}
\item\label{prof1}
The action $\sigma$ is profinite.
\item\label{prof2}
The action $\sigma$ is compact and the von Neumann algebra $P$ is
AFD.
\item\label{prof3}
The action $\sigma$ is weakly compact.
\item\label{prof4}
There exists a state $\p$ on $\B(L^2(P))$ such that
$\p|_P=\tau$ and $\p\circ\Ad u=\p$ for all $u\in\U(P)\cup\sigma(\G)$.
\item\label{prof5}
The von Neumann algebra $L(\G)$ is co-amenable in $P\rtimes\G$.
\end{enumerate}
Then, one has
$(\ref{prof1})\Rightarrow(\ref{prof2})\Rightarrow(\ref{prof3})
\Leftrightarrow(\ref{prof4})\Leftrightarrow(\ref{prof5})$.
\end{prop}
(Note that, by a result of H{\o}egh-Krohn--Landstad--St{\o}rmer (\cite{hls}),
if in the above statement we restrict our attention to ergodic actions
$\Gamma \curvearrowright P$, then the condition that $P$ is AFD in
part $(2)$ follows automatically from
the assumption $\Gamma\curvearrowright P$ compact.
We observe that weak compactness also implies that $P$ is AFD by
Connes's theorem (\cite{connes:cls}).)
\begin{proof}
We have $(\ref{prof1})\Rightarrow(\ref{prof2})$, by the
definitions. We prove $(\ref{prof2})\Rightarrow(\ref{prof4})$.
Since $P$ is AFD, there is a net $\Phi_n$ of normal u.c.p.\
maps from $\B(L^2(P))$ into $P$ such that
$\tau\circ(\Phi_n|_P)=\tau$ and $\|a-\Phi_n(a)\|_2\to0$
for all $a\in P$.
Let $G$ be the SOT-closure of $\sigma(\G)$ in the unitary group on $L^2(P)$.
By assumption, $G$ is a compact group and has a normalized
Haar measure $m$. We define a state $\p_n$ on $\B(L^2(P))$ by
\[
\p_n(x)=\int_G \tau\circ\Phi_n(gxg^{-1})\,dm(g).
\]
It is clear that $\p_n$ is $\Ad(\G)$-invariant and $\p_n|_P=\tau$.
We will prove that the net $\p_n$ is approximately $P$-central.
Let $\Phi_n(x)=V^*\pi(x)V$ be a Stinespring dilation.
Then, for $x\in\B(L^2(P))$ and $a\in P$, one has
\begin{align*}
\|\Phi_n(xa)-\Phi_n(x)\Phi_n(a)\|_2
&=\|V^*\pi(x)(1-VV^*)\pi(a)V\widehat{1}\|_{L^2(P)}\\
&\le\|x\|\|(1-VV^*)^{1/2}\pi(a)V\widehat{1}\|_{L^2(P)}\\
&=\|x\|\tau(\Phi_n(a^*a)-\Phi_n(a^*)\Phi_n(a))^{1/2}\\
&\le2\|x\|\|a\|^{1/2}\|a-\Phi_n(a)\|_2^{1/2}.
\end{align*}
It follows that for every $x\in\B(L^2(P))$ and $a\in P$, one has
\[
|\p_n(xa)-\p_n(ax)|\le 4\|x\|\|a\|^{1/2}\sup_{g\in G}\|gag^{-1}-\Phi_n(gag^{-1})\|_2^{1/2},
\]
which converge to zero since $\{ gag^{-1} : g\in G\}$ is compact in $L^2(P)$
and $\Phi_n$'s are contractive on $L^2(P)$.
Hence $\p_n$ is approximately $P$-central and $\p=\Lim_n\p_n$ satisfies the requirement.

We prove $(\ref{prof3})\Leftrightarrow(\ref{prof4})$.
Take a net $\eta_n$ satisfying the conditions $(1)$, $(2)$ and $(3')$
of Definition \ref{defn:profinite}.
We define a state $\p$ on $\B(L^2(P))$ by $\p=\Lim_n\p_n$ with
$\p_n(x)=\ip{(x\otimes1)\eta_n,\eta_n}$.
Then, for any $u\in\U(P)\cup\sigma(\G)$, one has
\[
\p(u^*xu)=\Lim_n\ip{(x\otimes1)(u\otimes\bar{u})\eta_n,(u\otimes\bar{u})\eta_n}
=\p(x)
\]
by the conditions $(1)$ and $(2)$ of Definition \ref{defn:profinite}.
That $\p|_P=\tau$ follows from $(3')$.
Conversely, suppose now that $\p$ is given.
We recall that $\B(L^2(P))$ is canonically identified with
the dual Banach space of the space $S_1(L^2(P))$ of trace
class operators. Take a net of positive elements $T_n\in S_1(L^2(P))$
with $\Tr(T_n)=1$ such that $\Tr(T_nx)\to\p(x)$ for every $x\in\B(L^2(P))$.
Let $b_n\in L^1(P)_+$ be such that $\Tr(T_na)=\tau(b_na)$
for $a\in P$. Since $\Tr(T_na)\to\p(a)=\tau(a)$ for $a\in P$,
the net $(b_n)$ converges to $1$ weakly in $L^1(P)$.
Thus, by the Hahn-Banach separation theorem, one may assume,
by passing to a convex combinations, that $\|b_n-1\|_1\to0$.
By a routine perturbation argument, we may assume further that $b_n=1$.
For the reader's convenience we give an argument for this.
Let $h(t)=\max\{1,t\}$ and $k(t)=\max\{1-t,0\}$ be functions on $[0,\infty)$,
and let $c_n=h(b_n)^{-1}$. We note that $0\le c_n\le1$ and $b_nc_n+k(b_n)=1$.
We define $T_n'=c_n^{1/2}T_nc_n^{1/2}+k(b_n)^{1/2}P_0k(b_n)^{1/2}$,
where $P_0$ is the orthogonal projection onto $\C\widehat{1}$.
Then, one has
\begin{align*}
\|T_n-T_n'\|_1 &\le2\|T_n^{1/2}-c_n^{1/2}T_n^{1/2}\|_2+\|k(b_n)\|_1\\
&=2\tau(b_n(1-c_n^{1/2})^2)^{1/2}+\|k(b_n)\|_1\\
&\le2\tau(b_n(1-c_n))^{1/2}+\|k(b_n)\|_1\\
&\le2\|b_n-1\|_1^{1/2}+\|1-b_n\|_1\to0.
\end{align*}
Hence, by replacing $T_n$ with $T_n'$,
we may assume that $\Tr(T_na)=\tau(a)$ for $a\in P$.
Since for every $x\in\B(L^2(P))$ and $u\in\U(P)\cup\sigma(\G)$,
one has
\[
\Tr((T_n-\Ad(u)T_n)x)\to\p(x)-\p(\Ad(u^*)(x))=0,
\]
by applying the Hahn-Banach separation theorem again,
one may furthermore assume that $\|T_n-\Ad(u)(T_n)\|_{S_1}\to0$
for every $u\in\U(P)\cup\sigma(\G)$.
Then by (\ref{eq:ps}), the Hilbert-Schmidt operators
$T_n^{1/2}$ satisfy
$\|T_n^{1/2}-\Ad(u)(T_n^{1/2})\|_{S_2}\to0$
for every $u\in\U(P)\cup\sigma(\G)$.
Now, if we use the standard identification between
$S_2(L^2(P))$ and $L^2(P\vt\bar{P})$ given by
\[
S_2(L^2(P))\ni \sum_k\ip{\,\cdot\,,\eta_k}\xi_k\mapsto
\sum_k\xi_k\otimes\bar{\eta}_k\in L^2(P\vt\bar{P})
\]
and view $T_n^{1/2}$ as an element $\zeta_n\in L^2(P\vt\bar{P})$,
then we have $\ip{(a\otimes1)\zeta_n,\zeta_n}=\tau(a)
=\ip{\zeta_n,(1\otimes\bar{a})\zeta_n}$ and
$\|\zeta_n-(u\otimes\bar{u})\zeta_n\|_2\to0$
for every $u\in\U(P)\cup\sigma(\G)$.
Therefore, the net of
$\eta_n=(\zeta_n\zeta_n^*)^{1/2}\in L^2(P\vt\bar{P})_+$
verifies the conditions of weak compactness.

Finally, we prove $(\ref{prof4})\Leftrightarrow(\ref{prof5})$.
We consider $P\rtimes\G$ as the von Neumann subalgebra of
$\B(L^2(P)\vt\ell^2(\G))$ generated by $P\vt\C1$ and
$(\sigma\otimes\lambda)(\G)$. This gives an identification between
$L^2(P\rtimes\G)$ and $L^2(P)\vt\ell^2(\G)$. Moreover, the basic
construction $\ip{P\rtimes\G,e_{L(\G)}}$ becomes $\B(L^2(P))\vt
L(\G)$, since it is the commutant of the right $L(\G)$-action
(which is given by $(1\otimes\rho)(\G)$).
Now suppose that $\p$ is given as in the condition $(\ref{prof4})$.
Then, $\tilde{\p}=\p\otimes\tau$ on $\B(L^2(P)\vt\ell^2(\G))$
is $\Ad(\U(P\vt\C1)\cup(\sigma\otimes\lambda)(\G))$-invariant
and $\tilde{\p}|_{P\rtimes\G}=\tau$.
This implies that $L(\G)$ is co-amenable in $P\rtimes\G$.
Conversely, if $\tilde{\p}$ is a $(P\rtimes\G)$-central state
such that $\tilde{\p}|_{P\rtimes\G}=\tau$, then
the restriction $\p$ of $\tilde{\p}$ to $\B(L^2(P))$
satisfies the condition $(\ref{prof4})$.
\end{proof}
Note that by part $(7)$ in Theorem \ref{Lambdacb}, if $\Lambda(\Gamma)=1$ and
$\Gamma \curvearrowright P$ is a profinite action then
$\Lambda(P\rtimes \Gamma)=1$. More generally we have the following.
(Compare this with \cite{jolissaint}.)

\begin{cor}\label{criterion}
Let $\G$ be weakly amenable and $\Gamma\curvearrowright P$
be a weakly compact action on an AFD von Neumann algebra.
Then, $P\rtimes\Gamma$ has the completely bounded approximation property
and $\Lambda_{\cb}(P\rtimes\Gamma)=\Lambda_{\cb}(\G)$.
\end{cor}
\begin{proof}
By Proposition~\ref{profinite}, $L(\G)$ is co-amenable in $P\rtimes\Gamma$.
Hence, Theorem 4.9 of \cite{delaroche} implies that
$\Lambda_{\cb}(P\rtimes\Gamma)=\Lambda_{\cb}(L(\G))=\Lambda_{\cb}(\G)$.
\end{proof}

\begin{prop}\label{weaklycompact} Let $P\subset M$
be an inclusion of finite von Neumann algebras such that $P'\cap M
\subset P$. Assume the normalizer $\mathcal N_M(P)$ contains a
subgroup $\mathcal G$ such that its action on $P$ is weakly compact
and $(P \cup \mathcal G)'' = \mathcal N_M(P)''$. Then the action of
$\mathcal N_M(P)$ on $P$ is weakly compact. Moreover, if $\mathcal
N_M(P)\curvearrowright P$ is weakly compact and $p\in \mathcal P(P)$
then $\mathcal N_{pMp}(pPp)\curvearrowright pPp$ is weakly compact.
\end{prop}
\begin{proof} We may clearly assume $\mathcal N_M(P)''=M$.
Denote by $\sigma$ the action of $\mathcal N_M(P)$ on $P$. If $u\in
\mathcal N_M(P)$, then by the conditions $P'\cap M=\mathcal Z(P)$
and $(P \cup \mathcal G)'' = M$ it follows that there exists a
partition $\{p_i\}_i \subset \mathcal \mathcal Z(P)$ and unitary
elements $v_i \in P$ such that $u=\Sigma_i p_iv_i u_i$ for some $v_i
\in \mathcal G$ (see e.g. \cite{dye}). Then $\sigma_v(x)=vxv^* =
\Sigma_i p_i \sigma_{v_iu_i}(x)$. Let now $\eta_n \in L^2(P \vt
P)_+$ satisfy the conditions in Definition 3.1 for the action
$\sigma_{|\mathcal G}$. By $3.1.(1)$ we have $\|\Sigma_i (p_i
\otimes p_i) \eta_n - \eta_n\|_2 \to0$, and thus $\|(p_i\otimes p_j)
\eta_n \|_2 \to0$, $\forall i\neq j$. Since
$q_i=\sigma_{u_i^*}(p_i)$ are mutually orthogonal as well, this also
implies that for $i\neq j$ we have
$$
\| (p_i \otimes
p_j)(\sigma_{v_iu_i}\otimes\bar{\sigma}_{v_ju_j})(\eta_n)\|_2
$$
$$
=\|(\sigma_{v_iu_i}\otimes\bar{\sigma}_{v_ju_j})((q_i \otimes q_j)
\eta_n)\|_2 =\|(q_i \otimes q_j) \eta_n\|_2 \to0
$$
Also, since $w_i=u_i^*v_iu_i \in \mathcal U(P)$, we have
$\|\sigma_{w_i}\otimes\bar{\sigma}_{w_i})(\eta_n)-\eta_n\|_2\to0$.
Combining with the condition $3.1.(2)$ on the action $\mathcal
G\curvearrowright P$, one gets $\|(p_i\otimes p_i) (\eta_n -
(\sigma_{v_iu_i}\otimes\bar{\sigma}_{v_iu_i})(\eta_n))\|_2 \to0$. By
Pythagora's theorem, all this entails
$$
\|\eta_n-(\sigma_v\otimes\bar{\sigma}_v)(\eta_n)\|^2_2  =
\Sigma_{i,j} \|(p_i\otimes p_j) \eta_n-(p_i \otimes p_j)
(\sigma_v\otimes\bar{\sigma}_v)(\eta_n)\|_2^2
$$
$$
=\Sigma_{i,j} \|(p_i\otimes p_j) \eta_n-(p_i \otimes p_j)
(\sigma_{v_iu_i}\otimes\bar{\sigma}_{v_ju_j})(\eta_n)\|_2^2 \to0,
$$
showing that $\mathcal N_M(P)\curvearrowright P$ satisfies
$3.1.(2)$, thus being weakly compact.

To see that weak compactness behaves well to reduction by
projections, note that any $v\in \mathcal N_{pMp}(pPp)$ extends to a
unitary in $\mathcal N_{M}(P)$. Thus, if $\p$ satisfies
the condition $\ref{profinite}.(\ref{prof4})$
for $\mathcal N_M(P)\curvearrowright P$ then
$\p^p=\p(p\,\cdot\,p)$ clearly satisfies the same condition
for $\mathcal N_{pMp}(pPp)\curvearrowright pPp$.
\end{proof}

The above result shows in particular that if a measure preserving
action of a countable group $\Gamma$ on a probability space
$(X,\mu)$ is weakly compact (i.e., $\Gamma \curvearrowright
L^\infty(X)$ weakly compact), then the action of its associated full
group $[\Gamma]$, as defined in \cite{dye}, is weakly compact. Thus,
weak compactness is an orbit equivalence invariant for group
actions, unlike profiniteness and compactness which are of course
not. In fact, by the last part of Proposition~\ref{weaklycompact},
it follows that weak compactness is a stably orbit equivalence (or
measure equivalence) invariant as well.

An embedding of finite von Neumann algebras $P\subset M$ is called
{\it weakly compact} if the action $\Nor_M(P) \curvearrowright
P$ is weakly compact. The next result shows that the complete metric
approximation property of a factor $M$ imposes the weak compactness
of all embeddings into $M$ of AFD (in particular abelian) von
Neumann algebras.

\begin{thm}\label{nssg}
Let $M$ be a finite von Neumann algebra with the c.m.a.p., i.e.\
$\Lambda_{\cb}(M)=1$. Then any embedding of an AFD von Neumann
algebra $P\subset M$ is weakly compact, i.e.,
$\Nor_M(P)\curvearrowright P$ is weakly compact, $\forall P\subset
M$ AFD subalgebra.
\end{thm}

For the proof, we need the following consequence of Connes's
Theorem \cite{connes:cls}. This is well-known, but we include
a proof for the reader's convenience.
\begin{lem}\label{extension}
Let $M$ be a finite von Neumann algebra,
$P\subset M$ be an AFD von Neumann subalgebra and $u\in\Nor_M(P)$.
Then, the von Neumann algebra $Q$ generated by $P$ and $u$ is AFD.
\end{lem}
\begin{proof}
Since $P$ is injective, the $\tau$-preserving conditional expectation
$E_P$ from $M$ onto $P$ extends to a u.c.p.\ map $\tilde{E}_P$
from $\B(L^2(M))$ onto $P$.
We note that $\tilde{E}_P$ is a conditional expectation:
$\tilde{E}_P(axb)=a\tilde{E}_P(x)b$ for every $a,b\in P$ and $x\in\B(L^2(M))$.
We define a state $\sigma$ on $\B(L^2(M))$ by
\[
\sigma(x)=\Lim_{n}
\frac{1}{n}\sum_{k=0}^{n-1}\tau(\tilde{E}_P(u^kxu^{-k})).
\]
It is not hard to check that $\sigma|_M=\tau$,
$\sigma\circ\Ad u=\sigma$ and  $\sigma\circ\Ad v=\sigma$ for every $v\in\U(P)$.
It follows that $\sigma$ is a $Q$-central state with $\sigma|_Q=\tau$.
By Connes's theorem, this implies that $Q$ is AFD.
\end{proof}

\begin{proof}[Proof of Theorem \ref{nssg}]
First we note the following general fact:
Let $\omega$ be a state on a C$^*$-algebra $N$ and $u\in\U(N)$.
We define $\omega_u(x)=\omega(xu^*)$ for $x\in N$.
Then, one has
\begin{equation}\label{eq:muomega0}
\max\{\|\omega-\omega_u\|,\|\omega-\omega\circ\Ad(u)\|\}\le2\sqrt{2|1-\omega(u)|}.
\end{equation}
Indeed, one has $\|\xi_\omega-u^*\xi_\omega\|^2=2(1-\Re\omega(u))\le2|1-\omega(u)|$,
where $\xi_\omega$ is the GNS-vector for $\omega$.

Let $(\phi_n)$ be a net of normal finite rank maps on $M$ such that
$\limsup\|\phi_n\|_{\cb}\le1$ and $\|x-\phi_n(x)\|_2\to0$ for all $x\in M$.
We observe that the net $(\tau\circ\phi_n)$ converges to $\tau$ weakly in $M_*$.
Hence by the Hahn-Banach separation theorem, one may assume,
by passing to convex combinations, that $\|\tau-\tau\circ\phi_n\|\to0$.
Let $\mu$ be the $*$-representation of the algebraic tensor product
$M\otimes\bar{M}$ on $L^2(M)$ defined by
\[
\mu(\sum_k a_k\otimes\bar{b}_k)\xi=\sum_k a_k\xi b_k^*.
\]
We define a linear functional $\mu_n$ on $M\otimes\bar{M}$ by
\[
\mu_n(\sum_k a_k\otimes\bar{b}_k)
=\ip{\mu(\sum_k \phi_n(a_k)\otimes\bar{b}_k)\widehat{1},\widehat{1}}_{L^2(M)}
=\tau(\sum_k \phi_n(a_k)b_k^*).
\]
Since $\phi_n$ is normal and of finite rank, $\mu_n$ extends to
a normal linear functional on $M\vt\bar{M}$, which is still denoted by $\mu_n$.
For an AFD von Neumann subalgebra $Q\subset M$, we denote by $\mu_n^Q$
the restriction of $\mu_n$ to $Q\vt\bar{Q}$.
Since $Q$ is AFD, the $*$-representation $\mu$ is continuous with respect to
the spatial tensor norm on $Q\otimes\bar{Q}$ and hence $\|\mu_n^Q\|\le\|\phi_n\|_{\cb}$.
We denote $\omega_n^Q=\|\mu_n^Q\|^{-1}|\mu_n^Q|$.
Since $\limsup\|\mu_n^Q\|\le1$ and $\lim\mu_n^Q(1\otimes1)=1$,
the inequality (\ref{eq:muomega0}), applied to $\omega_n^Q$, implies that
\begin{equation}\label{eq:muomega1}
\limsup_n\|\mu_n^Q-\omega_n^Q\|=0.
\end{equation}
Now, consider the case $Q=P$.
Since $\mu_n^P(v\otimes\bar{v})=\tau(\phi_n(v)v^*)\to1$ for any $v\in\U(P)$,
one has
\begin{equation}\label{eq:muomega2}
\limsup_n\|\omega_n^P-(\omega_n^P)_{v\otimes\bar{v}}\|=0
\end{equation}
by (\ref{eq:muomega0}) and (\ref{eq:muomega1}).
Now, let $u\in\Nor(P)$ and consider the case $Q=\ip{P,u}$,
which is AFD by Lemma~\ref{extension}.
Since $\mu_n^{\ip{P,u}}(u\otimes\bar{u})=\tau(\phi_n(u)u^*)\to1$,
one has
\begin{equation}\label{eq:muomega3}
\limsup_n\|\mu_n^{\ip{P,u}}-\mu_n^{\ip{P,u}}\circ\Ad(u\otimes\bar{u})\|=0
\end{equation}
by (\ref{eq:muomega0}) and (\ref{eq:muomega1}).
But since
$(\mu_n^{\ip{P,u}}\circ\Ad(u\otimes\bar{u}))|_{P\vt\bar{P}}=\mu_n^P\circ\Ad(u\otimes\bar{u})$,
one has
\begin{equation}\label{eq:muomega4}
\limsup_n\|\omega_n^P-\omega_n^P\circ\Ad(u\otimes\bar{u})\|=0
\end{equation}
by (\ref{eq:muomega1}) and (\ref{eq:muomega3}).
Now, we view $\omega_n^P$ as an $\zeta_n$ element in $L^1(P\vt\bar{P})_+$ and
let $\eta_n=\zeta_n^{1/2}$.
By (\ref{eq:ps}), the net $\eta_n$ satisfies all the required conditions.
\end{proof}
\section{Main Results}
We prove in this section the main results of the paper. They will
all follow from the following stronger version of
the Theorem stated in the introduction:

\begin{thm}\label{crossedprod}
Let $\G=\F_{r(1)}\times\cdots\times\F_{r(k)}$ be a direct product of
finitely many free groups of rank $2\le r(j)\le\infty$ and
denote by $\G_j$ the kernel of the projection from $\G$ onto $\F_{r(j)}$.
Let $M=Q\rtimes\G$ be the crossed product of
a finite von Neumann algebra $Q$ by $\G$ (action need not be ergodic nor free).
Let $P\subset M$ be such that $P\not\preceq_M Q$.
Let $\GG\subset\Nor_M(P)$ be a subgroup which acts weakly compactly
on $P$ by conjugation, and denote $N=\GG''$.
Then there exist projections $p_1,\ldots,p_k\in\Zt(N'\cap M)$
with $\sum_{j=1}^kp_j=1$ such that $Np_j\lessdot_MQ\rtimes\G_j$ for every $j$.
\end{thm}

>From the above result, we will easily deduce several
(in)decomposability properties for certain factors constructed out
of free groups and their profinite actions. Note that Corollaries
\ref{cor:free} and \ref{cor:freetensor} below are just Corollaries 1
and 2 in the introduction, while Corollary  \ref{cor:uniquecartan}
is a generalization of Corollary 3 therein.

\begin{cor}\label{cor:free}
If $P\subset L(\Bbb F_r)^t$ is a diffuse AFD von Neumann subalgebra
of the amplification by some $t>0$
of a free group factor $L(\Bbb F_r), 2\leq r \leq \infty$,
then $\Nor_{L(\F_r)}(P)''$ is AFD.
\end{cor}

Note that the above corollary generalizes the in-decomposability
results for free group factors in \cite{ozawa:solid} and
\cite{voiculescu}. Indeed, Voiculescu's celebrated result in
\cite{voiculescu}, showing that the normalizer of any amenable
diffuse subalgebra $P\subset L(\F_r)$ cannot generate all $L(\F_r)$,
follows from \ref{cor:free} because $L(\F_r)$ is non-AFD by \cite{mvn}. Also,
since any unitary element commuting with a subalgebra $P\subset N$
lies in the normalizer of $P$, \ref{cor:free} shows in particular that the
commutant of any diffuse AFD subalgebra $P\subset N=L(\F_r)$ is
amenable, i.e.\ $L(\F_r)$ is {\it solid} in the sense of
\cite{ozawa:solid}, which amounts to the free group case of a result
in \cite{ozawa:solid}. Note however that the in-decomposability results
in \cite{voiculescu} and \cite{ozawa:solid} cover much larger classes
of factors, e.g. all free products of diffuse von Neumann
algebras in \cite{voiculescu} (for absence of
Cartan subalgebras) and all II$_1$ factors arising from
word-hyperbolic groups in \cite{ozawa:solid} (for solidity).

Calling {\it strongly solid} (or {\it
s-solid}) the factors satisfying the property that the normalizer of
any diffuse amenable subalgebra generates an amenable von Neumann
algebra, one can at this point speculate that any c.m.a.p. s-solid
factor may in fact follow isomorphic to an amplification of a free
group factor (i.e., to an interpolated free group factor
\cite{dykema}, \cite{radulescu}).

\begin{cor}\label{cor:freetensor}
If $Q$ is a type $\mathrm{II}_1$-factor with c.m.a.p., then $Q\vt
L(\F_r)$ does not have Cartan subalgebras. Moreover, if $M \subset
Q\vt L(\F_r)$ is a subfactor of finite index, then $M$ does not have
Cartan subalgebras either.
\end{cor}

This corollary shows in particular that if $Q$ is an arbitrary
subfactor of a tensor product of free group factors, then $Q\vt
L(\F_r)$ (or any of its finite index subfactors) has no Cartan
subalgebras. When applied to $Q=R$, this shows that the subfactor $M
\subset R\vt L(\F_r)$ with $M\not\simeq M^{\op}$ constructed in
\cite{connes:chi}, as the fixed point algebra of an appropriate free
action of a finite group on $R\vt L(\F_r)$ (which thus has finite
index in $R\vt L(\F_r)$), does not have Cartan subalgebras.

Another class of factors without Cartan subalgebras is provided by
part $(2)$ of the next corollary.

\begin{cor}\label{cor:sr}
Let $\Gamma = \F_{r(1)}\times\cdots\times\F_{r(k)}$, as in $\ref{crossedprod}$,
and $\Gamma \curvearrowright X$ an ergodic
probability-measure-preserving action. Then $M=L^\infty(X)\rtimes
\Gamma$ is a II$_1$ factor and for each $t>0$ we have:
\begin{enumerate}
\item\label{sr1}
Assume $M^t$ has a maximal abelian $^*$-subalgebra $A$
such that $\mathcal N_{M^t}(A) \curvearrowright A$ is weakly compact
and $N=\mathcal N_{M^t}(A)''$ is a subfactor of finite index in
$M^t$. Then $\Gamma \curvearrowright X$ is necessarily a free
action, $L^\infty(X)$ is Cartan in $M$ and there exists a unitary
element $u\in M^t$ such that $uAu^* = L^\infty(X)^t$.
\item\label{sr2}
Assume $\Gamma \curvearrowright X$ is profinite (or merely compact).
Then $M$ has a Cartan subalgebra if and only if $\Gamma
\curvearrowright X$ is free.
\item\label{sr3}
Assume $\Gamma=\Bbb F_r$. If $M^t$ has a weakly compact maximal
abelian $^*$-subalgebra $A$ whose normalizer generates a von Neumann
algebra without amenable direct summand. Then $\Gamma
\curvearrowright X$ follows free and $A$ is unitary conjugate to
$L^\infty(X)^t$.
\end{enumerate}
\end{cor}

Note that one can view part $(1)$ of the above corollary as a {\it
strong rigidity result}, in the spirit of results in
(\cite{popa:betti}, \cite{popa:strong}, \cite{ipp}). Indeed, by
taking $A=L^\infty(Y)$ to be Cartan in $M^t$, it follows that any
isomorphism between group measure space II$_1$ factors
$\theta:(L^\infty(X)\rtimes \Gamma)^t \simeq L^\infty(Y) \rtimes
\Lambda$, with the ``source'' $\Gamma$ a direct product of finitely
many free groups and the ``target'' $\Lambda$ arbitrary but the
action $\Lambda \curvearrowright Y$ weakly compact (e.g. profinite,
or compact), is implemented by a stable orbit equivalence of the
free ergodic actions $\Gamma \curvearrowright X$, $\Lambda
\curvearrowright Y$, up to perturbation by an inner automorphism and
by an automorphism coming from a 1-cocycle of the target action.

\begin{cor}\label{cor:uniquecartan}
Let $\Gamma = \F_{r(1)}\times\cdots\times\F_{r(k)}$ (as in $\ref{crossedprod},
\ref{cor:sr}$) and $\Gamma \curvearrowright X$ a free ergodic profinite (or
merely compact) action. Then, $L^\infty(X)$ is the unique Cartan
subalgebra of the $\mathrm{II}_1$-factor $L^\infty(X)\rtimes\Gamma$,
up to unitary conjugacy. Moreover, if $\mathcal F\mathcal P$ denotes
the class of all II$_1$ factors that can be embedded as subfactors
of finite index in some $L^\infty(X)\rtimes \Gamma$, with $\Gamma
\curvearrowright X$ free ergodic compact action and $\Gamma$ as
above, then any $M\in \mathcal F\mathcal P$ has unique Cartan
subalgebra, up to unitary conjugacy. The class $\mathcal F\mathcal
P$ is closed to amplifications, tensor product and finite index
extension/restriction. Also, if $M\in \mathcal F\mathcal P$ and
$N\subset M$ is an irreducible subfactor of finite index, then
$[M:N]$ is an integer.
\end{cor}

The above corollary implies that any isomorphism between factors
$M\in \mathcal F \mathcal P$ comes from an isomorphism of the orbit
equivalence relations $\mathcal R_M$ associated with their unique
Cartan decomposition. Hence, like in the case of the class $\mathcal
H\mathcal T$ of factors in \cite{popa:betti}, invariants of
equivalence relations, such as Gaboriau's cost and $L^2$-Betti
numbers (\cite{gaboriau}), are isomorphism invariants of II$_1$
factors in $\mathcal F\mathcal P$. Note that all factors in the
class $\mathcal F\mathcal P$ have $\Lambda_{\cb}$-constant equal to
1 by Theorem \ref{Lambdacb} and have Haagerup's compact
approximation property by \cite{haagerup:map}.

The sub-class of II$_1$ factors $L^\infty(X)\rtimes \Bbb F_r\in
\mathcal F\mathcal P$, arising from free ergodic profinite
probability-measure-preserving actions of free groups $\Bbb F_r
\curvearrowright X$, is of particular interest, as they are
inductive limits of (amplifications of) free group factors. We call
such a factor $L^\infty(X)\rtimes\F_r$ an \emph{approximate free
group factor} of rank $r$. By Corollary~\ref{cor:uniquecartan}, more
than being in the class $\mathcal F\mathcal T$, such a factor has
the property that any maximal abelian $^*$-subalgebra with
normalizer generating a von Neumann algebra with no amenable summand
is unitary conjugate to $L^\infty(X)$. When combined with
\cite{gaboriau}, we see that approximate free group factors of
different rank are not isomorphic and that for $r<\infty$ they have
trivial fundamental group. Also, they are prime by
\cite{ozawa:kurosh}, in fact by \ref{crossedprod} the normalizer (in particular
the commutant) of any AFD II$_1$ subalgebra of such a factor must
generate an AFD von Neumann algebra. We will construct uncountably
many approximate free group factors in Section \ref{sec:uncountable}
and comment more on this class in Remark \ref{rem:appfree}.

For the proof of Theorem~\ref{crossedprod}, recall from
\cite{popa:bernoulli,popa:indec} the construction of $1$-parameter
automorphisms $\alpha_t$ (``malleable deformation'') of
$L(\F_r\ast\widetilde{\F}_r)$. Let $\widetilde{\F}_r$ be a copy of $\F_r$
and $a_1,a_2,\ldots$ (resp.\ $b_1,b_2,\ldots$) be
the standard generators of $\F_r$ (resp.\ $\widetilde{\F}_r$) viewed
as unitary elements in $L(\F_r\ast\widetilde{\F}_r)$.
Let $h_s=(\pi\sqrt{-1})^{-1}\log b_s$, where $\log$ is the principal
branch of the complex logarithm so that $h_s$ is a self-adjoint
element with spectrum contained in $[-1,1]$.
For simplicity, we write $b_s^t$ ($s=1,2,\ldots$ and $t\in\R$)
for the unitary element $\exp(t\pi\sqrt{-1}h_s)$.
The $*$-automorphism $\alpha_t$ is defined by $\alpha_t(a_s)=b_s^ta_s$
and $\alpha_t(b_s)=b_s$.

In this paper, we adapt this construction
to $\G=\F_{r(1)}\times\cdots\times\F_{r(k)}$ acting on $Q$ and $M=Q\rtimes\G$.
We extend the action $\G\curvearrowright Q$ to that of
\[
\widetilde{\G}=(\F_{r(1)}\ast\widetilde{\F}_{r(1)})
\times\cdots\times(\F_{r(k)}\ast\widetilde{\F}_{r(k)}),
\]
where $\widetilde{\F}_{r(j)}$'s act trivially on $Q$.
We denote by $a_{j,1},a_{j,2},\ldots$ (resp.\ $b_{j,1},b_{j,2},\ldots$)
the standard generators of $\F_{r(j)}$ (resp.\ $\widetilde{\F}_{r(j)}$)
We redefine the $*$-homomorphism
\[
\alpha_t\colon M\to\widetilde{M}=Q\rtimes\widetilde{\G}
\]
by $\alpha_t(x)=x$ for $x\in Q$ and
$\alpha_t(a_{j,s})=b_{j,s}^ta_{j,s}$ for each $1\le j\le k$ and $s$.
(We can define $\alpha_t$ on $\widetilde{M}$, but we do not need it.)

Let
\[
\gamma(t)=\tau(b_{j,s}^t)
=\frac{1}{2}\int_{-1}^1\exp(t\pi\sqrt{-1}h)\,dh
=\frac{\sin(t\pi)}{t\pi}=\gamma(-t)
\]
and $\phi_{j,\gamma(t)}\colon L(\F_{r(j)})\to L(\F_{r(j)})$ be
the Haagerup multiplier (\cite{haagerup:map}) associated
with the positive type function $g\mapsto\gamma(t)^{|g|}$ on $\F_{r(j)}$.
We may extend
\[
\phi_{\gamma(t)}=\phi_{1,\gamma(t)}\otimes\cdots\otimes\phi_{k,\gamma(t)}
\]
to $M$ by defining
$\phi_{\gamma(t)}(x\lambda(g))=x\phi_{\gamma(t)}(\lambda(g))$
for $x\in Q$ and $\lambda(g)\in L(\G)$.
We relate $\alpha_t$ and $\phi_{\gamma(t)}$ as follows (cf.\ \cite{peterson}).

\begin{lem}\label{haagerupmult}
One has $E_M\circ\alpha_t=\phi_{\gamma(t)}$.
\end{lem}

\begin{proof}
Since $E_M(x\lambda(g))=xE_{L(\G)}(\lambda(g))$
for $x\in Q$ and $\lambda(g)\in L(\widetilde{\G})$,
one has $E_M\circ\alpha_t(x\lambda(g))=xE_{L(\G)}(\alpha_t(\lambda(g)))$ for
$x\in Q$ and $\lambda(g)\in L(\G)$.
Hence it suffices to show $E_{L(\G)}\circ\alpha_t=\phi_{\gamma(t)}$
on $L(\G)$.
Since all $E_{L(\G)}$, $\alpha_t$ and $\phi_{\gamma(t)}$ split as tensor products,
we may assume that $k=1$.
Since $a_1,\ldots,b_1,\ldots$ are mutually free,
it is not hard to check
\[
(E_{L(\F_r)}\circ\alpha_t)(a_{i_1}^{\pm1}\cdots a_{i_n}^{\pm1})
=\gamma(t)^n a_{i_1}^{\pm1}\cdots a_{i_n}^{\pm1}
=\phi_{\gamma(t)}(a_{i_1}^{\pm1}\cdots a_{i_n}^{\pm1})
\]
for every reduced word $a_{i_1}^{\pm1}\cdots a_{i_n}^{\pm1}$ in $\F_r$.
\end{proof}

In particular, the u.c.p.\ map $E_M\circ\alpha_t$ on $M$
is compact over $Q$ provided that $r(j)<\infty$ for every $j$.
In case of $r(j)=\infty$, we need a little modification:
we replace the defining equation $\alpha_t(a_{j,s})=b_{j,s}^ta_{j,s}$
with $\alpha_t(a_{j,s})=b_{j,s}^{st}a_{j,s}$.
Then, the u.c.p.\ map $E_M\circ\alpha_t$ is
compact over $Q$ and $\alpha_t\to\id_M$ as $t\to0$.

Let $\G_j$ be the kernel of the projection from $\G$ onto $\F_{r(j)}$
and $Q_j=Q\rtimes\G_j\subset M$.
We consider the basic construction $\ip{M,e_{Q_j}}$ of $(Q_j\subset M)$.
Then, $L^2\ip{M,e_{Q_j}}$ is naturally an $M$-bimodule.

\begin{lem}\label{moduledecomp}
Let $Q_j\subset M\subset\widetilde{M}$ be as above.
Then, $L^2(\widetilde{M})\ominus L^2(M)$ is isomorphic as an $M$-bimodule
to a submodule of a multiple of $\bigoplus_{j=1}^k L^2\ip{M,e_{Q_j}}$.
\end{lem}
\begin{proof}
Let $\widetilde{\G}_j$ be the kernel of the projection from $\widetilde{\G}$
onto $\F_{r(j)}\ast\widetilde{\F}_{r(j)}$.
By permuting the position appropriately,
we consider that $\widetilde{\G}_j\times\F_{r(j)}\subset\widetilde{\G}$
and $\bigcap\widetilde{\G}_j\times\F_{r(j)}=\G$.
Let $\widetilde{Q}_j=Q\rtimes\widetilde{\G}_j$ and
$\widetilde{M}_j=Q\rtimes(\widetilde{\G}_j\times\F_{r(j)})$.
Since $L^2(M)=\bigcap_{j=1}^k L^2(\widetilde{M}_j)$,
it suffices to show $L^2(\widetilde{M})\ominus L^2(\widetilde{M}_j)$
is isomorphic as an $M$-bimodule to a multiple of $L^2\ip{M,e_{Q_j}}$.

We observe that
\[
L^2(\widetilde{M})\ominus L^2(\widetilde{M}_j)
=\bigoplus_d \Bigl[\widetilde{Q}_j\lambda(\F_{r(j)}d\,\F_{r(j)})\Bigr]
\]
where the square bracket means the $L^2$-closure and the
direct sum runs all over $d\in\F_{r(j)}\ast\widetilde{\F}_{r(j)}$
whose initial and final letters in the reduced form come from $\widetilde{\F}_{r(j)}$.
Let $\pi_j\colon\F_{r(j)}\ast\widetilde{\F}_{r(j)}\to\F_{r(j)}$ be the
projection sending $\widetilde{\F}_{r(j)}$ to $\{1\}$.
It is not difficult to see that
\[
x\lambda(gdh)\mapsto x\lambda(g)e_{Q_j}\lambda(\pi_j(d)h)
\]
extends to an $M$-bimodule isometry from
$\Bigl[\widetilde{Q}_j\lambda(\F_{r(j)}d\,\F_{r(j)})\Bigr]$
onto $L^2\ip{M,e_{Q_j}}$.
\end{proof}

We summarize the above two lemmas as follows.
\begin{prop}\label{prop:abstract}
Let $Q\subset Q_j\subset M$ be as above. Then,
there are a finite von Neumann algebra $\widetilde{M}\supset M$ and
trace-preserving $*$-homomorphisms $\alpha_t\colon M\to\widetilde{M}$
such that
\begin{enumerate}
\item $\lim_{t\to 0}\|\alpha_t(x)-x\|_2\to0$ for every $x\in M$;
\item $E_M\circ\alpha_t$ is compact over $Q$ for every $t>0$; and
\item $L^2(\widetilde{M})\ominus L^2(M)$ is isomorphic as an $M$-bimodule
to a submodule of a multiple of $\bigoplus_{j=1}^k L^2\ip{N,e_{Q_j}}$.
\end{enumerate}
\end{prop}

We complete the proof of Theorem~\ref{crossedprod} in this abstract setting.

\begin{thm}
Let $Q\subset Q_j\subset M$ be as in Proposition~$\ref{prop:abstract}$.
Let $P\subset M$ be such that $P\not\preceq_M Q$.
Let $\GG\subset\Nor_M(P)$ be subgroup which acts weakly compactly
on $P$ by conjugation, and $N=\GG''$.
Then there exist projection $p_1,\ldots,p_k\in\Zt(N'\cap M)$
with $\sum_{j=1}^kp_j=1$ such that $Np_j\lessdot_MQ_j$ for every $j$.
\end{thm}

\begin{proof}
We may assume that $\U(P)\subset\GG$.
We use Corollary~\ref{cor:relative} to conclude the relative amenability.
Let a non-zero projection $p$ in $\Zt(N'\cap M)$,
a finite subset $F\subset\GG$ and $\e>0$ be given arbitrary.
It suffices to find $\xi\in\bigoplus\bigoplus_{j=1}^k L^2\ip{N,e_{Q_j}}$
such that $\|x\xi\|_2\le\|x\|_2$ for all $x\in N$,
$\|p\xi\|_2\geq\|p\|_2/8$ and $\|[\xi,u]\|_2^2<\e$ for every $u\in F$.

Let $\delta=\|p\|_2/8$.
We choose and fix $t>0$ such that $\alpha=\alpha_t$
satisfies $\|p-\alpha(p)\|_2<\delta$ and
$\|u-\alpha(u)\|_2<\e/6$ for every $u\in F$.
We still denote by $\alpha$ when it is viewed as an isometry from
$L^2(M)$ into $L^2(\widetilde{M})$.
Let $(\eta_n)$ be the net of unit vectors in $L^2(P\vt \bar{P})_+$
as in Definition \ref{defn:profinite} and denote
\[
\widetilde{\eta}_n=(\alpha\otimes1)(\eta_n)\in L^2(\widetilde{M})\vt L^2(\bar{M}).
\]
We note that
\begin{equation}\label{eq:projective}
\|(x\otimes 1)\widetilde{\eta}_n\|_2^2
=\tau(\alpha^{-1}(E_{\alpha(M)}(x^*x)))=\|x\|_2^2
\end{equation}
for every $x\in\widetilde{M}$.
In particular, one has
\begin{equation}\label{eq:uualphaeta}
\|[u\otimes\bar{u},\widetilde{\eta}_n]\|_2
\le\|[u\otimes\bar{u},\eta_n]\|_2+2\|u-\alpha(u)\|_2
<\e/2
\end{equation}
for every $u\in F$ and large enough $n\in\N$.
We denote $\zeta_n=(e_{M}\otimes1)(\widetilde{\eta}_n)$ and
$\zeta_n^\perp=\widetilde{\eta}_n-\zeta_n$.
Noticing that $L^2(M)\vt L^2(\bar{M})$ is an $M\vt\bar{M}$-bimodule,
it follows form (\ref{eq:uualphaeta}) that
\begin{equation}\label{eq:uuzetan}
\|[u\otimes\bar{u},\zeta_n]\|_2^2
+\|[u\otimes\bar{u},\zeta_n^\perp]\|_2^2
=\|[u\otimes\bar{u},\widetilde{\eta}_n]\|_2^2<(\e/2)^2
\end{equation}
for every $u\in F$ and large enough $n\in\N$.
We claim that
\begin{equation}\label{eq:pzetaperp}
\Lim_n\|(p\otimes1)\zeta_n^\perp\|_2\geq\delta.
\end{equation}
Suppose this is not the case.
Then, for any $v\in\U(P)$, one has
\begin{align*}
&\Lim_n \|(p\otimes1)\widetilde{\eta}_n-(e_M\alpha(v)p\otimes\bar{v})\zeta_n\|_2\\
&\quad\le \Lim_n\|(p\otimes1)\widetilde{\eta}_n-(e_M\alpha(v)p\otimes\bar{v})
 \widetilde{\eta}_n\|_2+\Lim_n\|(p\otimes1)\zeta_n^\perp\|_2\\
&\quad\le \Lim_n\|(p\otimes1)\widetilde{\eta}_n-(e_Mp\otimes1)(\alpha(v)\otimes\bar{v})
 \widetilde{\eta}_n\|_2+\|[\alpha(v),p]\|_2+\delta\\
&\quad\le \Lim_n\|(p\otimes1)\zeta_n^\perp\|_2
 +\Lim_n\|\widetilde{\eta}_n-(\alpha(v)\otimes\bar{v})\widetilde{\eta}_n\|_2
 +2\|p-\alpha(p)\|_2+\delta\\
&\quad\le 4\delta
\end{align*}
since $pe_M=e_Mp$.
It follows that
\begin{equation}\label{eq:uniflb}\begin{aligned}
\|(E_M\circ\alpha)(v)p\|_2
&= \Lim_n \|\bigl((E_M\circ\alpha)(v)p \otimes\bar{v}\bigr)\widetilde{\eta}_n\|\\
&\geq \Lim_n \|(e_M\otimes1)\bigl((E_M\circ\alpha)(v)p
\otimes\bar{v}\bigr)\widetilde{\eta}_n\|\\
&= \Lim_{n}\|\bigl(e_M\alpha(v)p \otimes\bar{v}\bigr)\zeta_n\|\\
&\geq\|p\|_2-4\delta>0
\end{aligned}\end{equation}
for all $v\in\U(P)$. (One has
$\|(E_M\circ\alpha)(vp)\|_2\geq\|p\|_2-6\delta$ as well.) Since
$E_M\circ\alpha$ is compact over $Q$, this implies $P\preceq_MQ$ by
Theorem~\ref{cor:cptq}, contradicting the assumption. Thus by
(\ref{eq:uuzetan}) and (\ref{eq:pzetaperp}), there exists $n\in\N$
such that $\zeta=\zeta_n^\perp\in (L^2(\widetilde{M})\ominus
L^2(M))\vt L^2(\bar{M})$ satisfies
$|[u\otimes\bar{u},\zeta]\|_2<\e/2$ for every $u\in F$ and
$\|(p\otimes1)\zeta\|_2\geq\delta$. We note that for all $x\in M$,
\begin{equation}\label{eq:zetatrace}\begin{aligned}
\|(x\otimes1)\zeta\|_2^2
 =\|(e_M^\perp\otimes1)(x\otimes1)\widetilde{\eta}_n\|_2^2
 \le\|(x\otimes1)\widetilde{\eta}_n\|_2^2=\|x\|_2^2.
\end{aligned}\end{equation}
By Proposition~\ref{prop:abstract}, we may view $\zeta$
as a vector $(\zeta_i)$ in $\bigoplus_i L^2\ip{N,e_{Q_{j(i)}}}\vt L^2(\bar{M})$.
We consider $\zeta_i\zeta_i^*\in L^1(\ip{N,e_{Q_{j(i)}}}\vt\bar{M})$ and
define $\xi_i=((\id\otimes\tau)(\zeta_i\zeta_i^*))^{1/2}$
and then $\xi=(\xi_i)\in\bigoplus_i L^2\ip{N,e_{Q_{j(i)}}}$.
Then, the inequality (\ref{eq:zetatrace}) implies
\[
\|x\xi\|_2^2=\sum_i\tau(x^*x(\id\otimes\tau)(\zeta_i\zeta_i^*))
=\|(x\otimes1)\zeta\|_2^2\le\|x\|_2^2,
\]
and for all $x\in N$. In particular,
\[
\|p\xi\|_2=\|(p\otimes1)\zeta\|_2\geq\delta.
\]
Finally, by (\ref{eq:ps}), one has
\begin{equation*}\begin{alignedat}{2}
\|[\xi,u]\|_2^2 &= \sum_i\|\xi_i-(\Ad u)(\xi_i)\|_2^2
& &\le \sum_i\|\xi_i^2-(\Ad u)(\xi_i^2)\|_1\\
&\le \sum_i\|\zeta_i\zeta_i^*-\Ad(u\otimes\bar{u})(\zeta_i\zeta_i^*)\|_1
& &\le \sum_i2\|\zeta_i\|_2\|[u\otimes\bar{u},\zeta_i]\|_2\\
&\le 2\|\zeta\|_2\|[u\otimes\bar{u},\zeta]\|_2
& &< \e
\end{alignedat}\end{equation*}
for every $u\in F$.
\end{proof}
Before proving the corollaries to
Theorem \ref{crossedprod}, we mention one more result
in the spirit of \ref{crossedprod}. Its proof is similar
to the above, but requires more involved technique
from \cite{ipp}.

\begin{thm}
Let $M=M_1\ast M_2$ be the free product of finite von Neumann
algebras and $P\subset M$ be a von Neumann subalgebra such that
$P\not\preceq_M M_i$ for $i=1,2$. If the action of
$\GG\subset\Nor_N(P)$ on $P$ is weakly compact, then $\GG''$ is AFD.
\end{thm}
\begin{proof}
We follow the proof of Theorem~\ref{crossedprod}, but
use instead the deformation $\alpha_t$ given in
Lemma 2.2.2 in \cite{ipp}.
Let a non-zero projection $p$ in the center of $\Zt(\G'\cap M)$,
a finite subset $F\subset\GG$ and $\e>0$ be given arbitrary.
Since $P\not\preceq_M M_i$ for $i=1,2$, one has
\[
\lim_{t\to 0}\inf\{\|(E_M\circ\alpha_t)(vp)\|_2 : v\in\U(P)\}<(999/1000)\|p\|_2
\]
by Proposition 3.4 and Theorem 4.3 in \cite{ipp}. (N.B. This is
because Proposition 3.4 is the only part where the rigidity
assumption in Theorem 4.3 of that paper is being used.) Hence, if we
choose $\delta>0$ small enough and $t>0$ accordingly, then one
obtains as in the proof of Theorem~\ref{crossedprod} that
\[
\Lim_{n}\|(p\otimes1)\zeta_n^\perp\|_2\geq\delta
\]
for $\zeta_n^\perp=((1-e_M)\otimes1)\widetilde{\eta}_n
\in L^2(\widetilde{M}\ominus M)\vt L^2(\bar{M})$.
Since $L^2(\widetilde{M}\ominus M)$ is a multiple of $L^2(M\vt M)$
as an $M$-bimodule, one obtains
$\xi\in\bigoplus L^2(M\vt M)$ such that
$\|x\xi\|_2=\|\xi x\|_2\le\|x\|_2$ for all $x\in M$,
$\|p\xi\|_2\geq\delta$ and $\|[u,\xi]\|_2<\e$ for every $u\in F$.
This proves that $\GG''$ is AFD.
\end{proof}

\begin{proof}[Proof of Corollary 4.2.]
This is a trivial consequence of Theorems~\ref{nssg}
and \ref{crossedprod}.
\end{proof}

\begin{proof}[Proof of Corollary 4.3.] Suppose there is a
Cartan subalgebra $A\subset M$ where $M\subset N=Q\vt L(\F_r)$ is a
subfactor of finite index. Since $\F_r$ is non-amenable, $N$ is not
amenable relative to $Q$, so by Proposition 2.3, $M$ is not amenable
relative to $Q$ inside $N$. Hence, by Theorems~\ref{nssg} and
\ref{crossedprod}, one has $A\preceq_NQ$. By Theorem 2.4, this
implies there exist projections $p\in A'\cap N$, $q\in Q$, an
abelian von Neumann subalgebra $A_0\subset qQq$ and a non-zero
partial isometry $v\in N$ such that $p_0=vv^* \in p(A'\cap N)p$,
$q_0=v^*v \in A_0'\cap qNq$ and $v^*(Ap_0)v=A_0q_0$. Since
$Q=L(\F_r)'\cap N$, by ``shrinking'' $q$ if necessary we may clearly
assume $q=\bigvee \{uq_0u^* : u\in \U(L(\F_r))\}$. Since $L(\F_r)q$
is contained in  $(A_0q)'\cap qNq$, this implies $q_0$ has central
support $1$ in the von Neumann algebra $(A_0q)'\cap qNq$. But
$(A_0q_0)'\cap q_0Nq_0= v^*(A'\cap N)v$ by spatiality and since
$M\subset N$ has finite index, $A\subset A'\cap N$ has finite index
as well (in the sense of \cite{pimsner-popa}) so $A'\cap N$ is type
I, implying $(A_0q_0)'\cap q_0Nq_0$ type I, and thus $(A_0q)'\cap
qNq$ type I as well. But $L(\F_r)\simeq L(\F_r) q \subset
(A_0q)'\cap qNq$, contradiction.
\end{proof}

For the proof of Corollary 4.4, we'll need the following general
observation.

\begin{lem}\label{freenesscrit} Let $\Gamma$ be an ICC group
and $\Gamma \curvearrowright X$ an ergodic measure-preserving
action. Let $M= L^\infty(X) \rtimes \Gamma$. Then $M$ is a factor.
Moreover, $L^\infty(X)$ is maximal abelian (thus Cartan) in $M$ if
and only if there is a maximal abelian $^*$-subalgebra $A\subset M$
such that $A \preceq_M L^\infty(X)$.
\end{lem}

\begin{proof} The first part is well known, its proof being
identical to the Murray-von Neumann classical argument in
\cite{mvn}, showing that if a group $\Gamma$ is ICC then its group
von Neumann algebra $L(\Gamma)$ is a factor. For the second part,
denote $B=L^\infty(X)$ and let $A\subset M$ be maximal abelian
satisfying $A \preceq_M B$. Then there exists a non-zero partial
isometry $v\in M$, projections $p\in A=A'\cap M$, $q\in B$ and a
unital isomorphism $\theta$ of $Ap$ onto a unital subalgebra $B_0$
of $Bq$ such that $va = \theta(a)v$, $\forall a\in Ap$. Denoting
$q'=vv^* \in B_0'\cap qMq$, it follows that $q'(B_0'\cap
qMq)q'=(B_0q')'\cap q'Mq'$. Since by spatiality $B_0q'=vAv^*$ is
maximal abelian, this implies $q'(B_0'\cap qMq)q'=vAv^*$. Thus,
$B_0'\cap qMq$ has a type I direct summand. Since $(Bq)' \cap qMq$
is a subalgebra of $B_0'\cap pMp$, it follows that $B'\cap M$ has a
type I summand. Since $\Gamma$ acts ergodically on $\mathcal
Z(B'\cap M)\supset B$ (or else $M$ wouldn't be a factor), the
algebra $B'\cap M$ is homogeneous of type I$_n$, for some
$n<\infty$.

Note at this point that since all maximal abelian subalgebras of the
type I summand of $B_0'\cap qMq$ containing $q'$ are unitary
conjugate (cf. \cite{kadison}), we may assume $q'$ is in a maximal
abelian algebra containing $Bq$. Thus, if $\mathcal Z$ denotes the
center of $B'\cap M$, then $\mathcal Zq'\subset q'(B_0'\cap
qMq)q'=B_0q'\subset Bq'$, showing that $\mathcal Zq'=Bq'$. Since $B,
\mathcal Z$ are $\Gamma$-invariant with the corresponding
$\Gamma$-actions ergodic, it follows that there exists a partition
of 1 with projections of equal trace $p_1, ... ,p_m\in \mathcal Z$
such that $\mathcal Z=\Sigma_i Bp_i$ and $E_B(p_i)=m^{-1}1$,
$\forall i$. Since $B'\cap M=\mathcal Z'\cap M$ has an orthonormal
basis over $\mathcal Z$ with $n^2$ unitary elements, this shows that
$B'\cap M$ has a finite unitary orthonormal basis over $B$. But if
$x\in (B'\cap M) \setminus B,$ and $x=\Sigma_g a_g u_g$ is its
Fourier series, with $a_g\neq 0$ for some $g\neq e$, then $p_gu_g\in
B'\cap M$, where $p_g$ denotes the support projection of $a_g$. Now,
since $\Gamma$ is ICC there exist infinitely many $h_n\in \Gamma$
such that $g_n=h_ngh_n^{-1}$ are distinct. This shows that all
$\sigma_{h_n}(p_g)u_{g_n}\subset B'\cap M$ are mutually orthogonal
relative to $B$. By \cite{pimsner-popa}, this contradicts the
finiteness of the index of $B\subset B'\cap M$. Thus, we must have
$B'\cap M=B$, showing that $\Gamma \curvearrowright X$ is free and
$B=L^\infty(X)$ is maximal abelian, hence Cartan.
\end{proof}

\begin{proof}[Proof of Corollary 4.4.] The factoriality of $M$ was shown in
\ref{freenesscrit} above.

To prove part $(1)$, note that $\mathcal N_{M^t}(A) \curvearrowright
A$ weakly compact implies $\mathcal N_{M}(A^{1/t}) \curvearrowright
A^{1/t}$ weakly compact, where $A^{1/t} \subset M^{1/t}$ is the
semiregular maximal abelian $^*$-subalgebra obtained by amplifying
$A\subset M$ by $1/t$. Since weak compactness behaves well to
amplifications (see comments after the proof of
\ref{weaklycompact}), this shows that it is sufficient to prove the
case $t=1$. Let $\G_j$ be as in \ref{crossedprod}. If $N=\mathcal
N_M(A)''\lessdot_M L^\infty(X) \rtimes\G_j$ for some $j$, then by
$[M:N]< \infty$ it follows that $M \lessdot_M L^\infty(X)
\rtimes\G_j$ as well. But this implies $\Bbb F_{r(j)}$ amenable, a
contradiction. Thus, by Theorem \ref{crossedprod} we have $A \preceq
L^\infty(X)$ and the statement follows from Lemma
\ref{freenesscrit}.

Part $(2)$ follows trivially from part $(1)$, since $\Gamma
\curvearrowright X$ compact implies $M$ has c.m.a.p., by
\ref{profinite}.

Arguing as in the proof of $(1)$, we see that to prove $(3)$ it is
sufficient to settle the case $t=1$. If $N \lessdot_M L^\infty(X)$,
then $N$ would follow amenable. Thus, by \ref{crossedprod} we have
$A \preceq L^\infty(X)$ and Lemma \ref{freenesscrit} above applies
again.
\end{proof}

The proof of Corollary 4.5 will follow readily from the next general
``principle''.

\begin{prop}\label{cartancrit} Assume a II$_1$ factor $M$ has the
property:

\vskip .05in \noindent $(a)$ $\exists$ $A\subset M$ Cartan and any
maximal abelian $^*$-subalgebra $A_0\subset M$ with $\mathcal
N_M(A_0)''$ a subfactor of finite index in $M$ is unitary conjugate
to $A$. \vskip .05in

Then any amplification and finite index extension/restriction of $M$
satisfies $(a)$ as well. Moreover, if $M$ satisfies $(a)$ and
$N\subset M$ is an irreducible subfactor of finite index, then
$[M:N]$ is an integer.
\end{prop}

\begin{proof} For the proof, we call an
abelian von Neumann subalgebra $B$ of a II$_1$ factor $P$ {\it
virtually Cartan} if it is maximal abelian and $Q=\mathcal N_P(B)''$
has finite dimensional center with $[qPq:Qq]<\infty$ for any atom
$q\in \mathcal Z(Q)$. We first prove that if $P\subset N$ is an
inclusion of factors with finite index and $B\subset P$ is virtually
Cartan in $P$ then any maximal abelian $^*$-subalgebra $A$ of
$B'\cap N$ is virtually Cartan in $N$.

To see this, note that that, by commuting squares, the index of
$B\subset B'\cap N$ (in the sense of \cite{pimsner-popa}) is
majorized by $[N:P]<\infty$, implying that $B'\cap N$ is a direct
sum of finitely many homogeneous type I$_{n_i}$ von Neumann algebras
$B_i$, with $1\leq n_1 < n_2 < ... < n_k<\infty$. Since any two
maximal abelian $^*$-subalgebras of a finite type I von Neumann
algebra are unitary conjugate and $\mathcal N_P(B)$ leaves $B'\cap
N$ globally invariant, it follows that given any $u\in N_P(B)$,
there exists $v(u)\in \mathcal U(B'\cap N)$ such that
$v(u)uAu^*v(u)^*=A$. Moreover, $A$ is Cartan in $B'\cap N$, i.e.
$\mathcal N_{B'\cap N}(A)''=B'\cap N$. This shows in particular that
the von Neumann algebra generated by $\mathcal N_N(A)$ contains
$B'\cap N$ and $v(u)u$, and thus it contains $u$, i.e. $\mathcal
N_P(B)\subset \mathcal N_N(A)''$. Thus, the
\cite{pimsner-popa}-index of $\mathcal N_N(A)''$ in $N$ is majorized
by the index of $P$ in $N$, and is thus finite. Since $N$ is a
factor, this implies $Q=\mathcal N_N(A)''$ has finite dimensional
center and $[qNq:Qq]<\infty$ for any atom in its center, i.e. $A$ is
virtually Cartan in $N$.

Now notice that since any unitary conjugacy of subalgebras $A,
A_0\subset M$ as in $(a)$ can be ``amplified'' to a unitary
conjugacy of $A^t, A_0^t$ in $M^t$, property $(a)$ is stable to
amplifications. This also shows that $(a)$ holds true for a factor
$M$ if and only if $M$ satisfies:

\vskip .1in \noindent $(b)$ {\it $\exists A \subset M$ Cartan and
any virtually Cartan subalgebra $A_0$ of $M$ is unitary conjugate to
$A$}. \vskip .1in

Since if a subfactor $N\subset M$ satisfies $[M:N]< \infty$ then
$\langle M, e_N \rangle$ is an amplification of $N$ (see e.g.
\cite{pimsner-popa}), it follows that in order to finish the proof
of the statement it is sufficient to prove that if $M$ satisfies
$(b)$ and $N\subset M$ is a subfactor with finite index, then $N$
satisfies $(b)$.

Let $A\subset M$ be a Cartan subalgebra of $M$. Let $P\subset N$ be
such that $N\subset M$ is the basic construction of $P\subset N$ (cf
\cite{jones}). Thus $P$ is isomorphic to an amplification of $M$ and
so it has a Cartan subalgebra $A_2\subset P$. By the first part of
the statement any maximal abelian subalgebra $A_1$ of $A_2'\cap N$
is virtually Cartan in $N$. Applying again the first part, any
maximal abelian $A_0$ of $A_1'\cap M$ is virtually Cartan in $M$, so
it is unitary conjugate to $A$. Thus, $A_0\subset M$ follows Cartan.
Thus, $L^2(M) = \oplus u_nL^2(A_0)$, for some partial isometries
$u_n\in M$ normalizing $A_0$. Since $A_0$ is a finitely generated
$A_1$-module, it follows that each $u_nL^2(A_0)$ is finitely
generated both as left and as right $A_1$ module, i.e. there exist
finitely many $\xi_i, \xi'_j \in u_nL^2(A_0)$ such that $\Sigma_i
\xi_i A_1$ and $\Sigma A_1 \xi'_i$ are dense in $u_nL^2(A_0)$. Thus,
if we denote by $\mathcal H_n$ the closure of the range of the
projection of $u_nL^2(A_0)$ onto $L^2(N)$ and by $\eta_i, \eta_j'$
the projection of $\xi_i, \xi_j$ onto $u_nL^2(A_0)$, then $\mathcal
H_n$ is a Hilbert $A_1$-bimodule generated as left Hilbert
$A_1$-module by $\eta_i\in L^2(N)$ and as a right Hilbert
$A_1$-module by $\eta_j'\in L^2(N)$. Moreover, since $\vee_n
u_nL^2(A_0)=L^2(M)$, we have $\vee_n \mathcal H_n = L^2(N)$. Thus,
by Section 1.4 in \cite{popa:betti}, $A_1$ is Cartan in $N$.

Note that the above argument shows that $N$ has Cartan subalgebra,
but also that any virtually Cartan subalgebra of $N$ is in fact
Cartan. If now $B_1\subset N$ is another Cartan subalgebra of $N$,
then let $B_0$ be a maximal abelian subalgebra of $B_1'\cap M$. By
the first part of the proof $B_0$ is virtually Cartan, so by $(b)$
there exists $v\in \mathcal U(M)$ such that $vA_0v^*=B_0$. Thus, if
we let $v_n=vu_n$ then $L^2(M)=\oplus_n v_n L^2(A_0)=\oplus_n
L^2(B_0)v_n$. Since $A_0$ (resp. $B_0$) is a finitely generated
$A_1$ (resp. $B_1$) module, there exist $\xi_i, \xi_j'\in
v_nL^2(A_0)=L^2(B_0)v_n$ such that $\Sigma_i \xi_i A_1$ is dense in
$v_nL^2(A_0)$ and $\Sigma_j B_1 \xi_j'$ is dense in $L^2(B_0)v_n$.
But then exactly the same argument as above shows that $L^2(N)$ is
spanned by Hilbert $B_1-A_1$ bimodules $\mathcal H_n$ which are
finitely generated both as right $A_1$ Hilbert modules and as left
Hilbert $B_1$ modules. By Section 1.4 in \cite{popa:betti}, it
follows that $A_1, B_1$ are unitary conjugate.

Finally, to see that for irreducible inclusions of factors $N\subset
M$ satisfying $(a)$ the index $[M:N]$ is an integer, when finite,
let $N\subset Q\subset P \subset M$ be the canonical intermediate
subfactors constructed in 7.1 of \cite{popa:betti}. Then $Q,P$
satisfy $(a)$ as well and by 7.1 in \cite{popa:betti} the Cartan
subalgebra of $P$ is maximal abelian and Cartan in $M$. Thus, as in
the proof of 7.2.3$^\circ$ in \cite{popa:betti}, we have $[Q:N],
[P:Q], [M:P]\in \Bbb N$, implying that $[M:N]\in \Bbb N$.
\end{proof}

\begin{proof}[Proof of Corollary 4.5.]
Let $M=L^\infty(X)\rtimes\Gamma$ and assume $A\subset M$ is a Cartan
subalgebra. By \ref{profinite} and \ref{criterion}, $M$ follows
c.m.a.p. Thus, \ref{nssg} applies to show that $\mathcal N_M(A)
\curvearrowright A$ is weakly compact. Since $\Bbb F_{r(j)}$ are all
non-amenable, $M=\mathcal N_M(A)''$ cannot be amenable relative to
$L^\infty(X) \rtimes \Gamma_j$ (with $\Gamma_j$ as defined in
\ref{crossedprod}), $\forall j$. Hence, Theorem \ref{crossedprod}
implies $A\preceq_ML^\infty(X)$. Then Lemma~\ref{conjugatecartan}
shows there is $u\in\U(M)$ such that $uAu^*=L^\infty(X)$, proving
the first part of the statement. The rest is
a consequence of Proposition \ref{cartancrit}.
\end{proof}

\section{Uncountably many approximate free group factors}\label{sec:uncountable}
In this section we prove that there are uncountably many approximate
free group factors of any rank $2 \leq n\leq \infty$. We do this by
using a ``separability argument,'' in the spirit of
\cite{popa:corresp,junge-pisier,ozawa:uncountable}. The proof is
independent of the previous sections. The result shows in particular
the existence of uncountably many orbit inequivalent profinite
actions of $\F_n$. The fact that $\F_n$ has uncountably many
orbit inequivalent actions was first shown in \cite{gaboriau-popa}.
A concrete family of orbit inequivalent actions of $\F_n$ was
recently obtained in \cite{ioana:uncountable}. Note that the actions
$\F_n \curvearrowright X$ in \cite{gaboriau-popa} and
\cite{ioana:uncountable} are not orbit equivalent to profinite
actions (because they have quotients that are free and have relative
property (T) in the sense of \cite{popa:betti}).

\begin{defn}
We say a unitary representation $(\pi,\hh)$ of $\G$
\emph{has (resp.\ essential) spectral gap} if there is a
finite subset $F$ of $\G$ and $\e>0$ such that the
self-adjoint operator
\[
\frac{1}{2|F|}\sum_{g\in F}(\pi(g)+\pi(g^{-1}))
\]
has (resp.\ essential) spectrum contained in $[-1,1-\e]$.
We say such $(F,\e)$ \emph{witnesses (resp.\ essential)
spectral gap} of $(\pi,\hh)$.
\end{defn}

It is well-known that $(\pi,\hh)$ has spectral gap
if and only if it does not contain approximate invariant
vectors.

\begin{defn}
Let $\G$ be a group. We say $\G$ is \emph{inner-amenable} (\cite{effros})
if the conjugation action of $\G$ on $\ell^2(\G\setminus\{1\})$
does not have spectral gap.

Let $\{\G_n\}$ be a family of finite index (normal) subgroups of $\G$.
We say $\G$ has the property $(\tau)$ with respect to $\{\G_n\}$ if
the unitary $\G$-representation on
\[
\bigoplus_n\ell^2(\G/\G_n)^{o}
\]
has spectral gap, where $\ell^2(\G/\G_n)^{o}=\ell^2(\G/\G_n)\ominus\C1_{\G/\G_n}$.

Let $I$ be a family of decreasing sequences
\[
i=\bigl(\G=\G^{(i)}_0\geq\G^{(i)}_1\geq\G^{(i)}_2\geq\cdots\bigr)
\]
of finite index normal subgroups of $\G$
such that $\bigcap\G^{(i)}_n=\{1\}$.
We allow the possibility that $\G^{(i)}_n=\G^{(i)}_{n+1}$.
We say the family $I$ is \emph{admissible} if
$\G$ has the property $(\tau)$ with respect to
$\{\G^{(i)}_m\cap \G^{(j)}_n : i,j\in I,\,m,n\in\N\}$
and
\[
\sup\{ [ \G : \G^{(i)}_m\G^{(j)}_n ] : m,n\in\N\}<\infty
\]
for any $i,j\in I$ with $i\neq j$.
\end{defn}

\begin{lem}
Let $\G\le\SL(d,\Z)$ with $d\geq2$ be a finite index subgroup and
\[
\G_n=\G\cap\ker\bigl(\SL(d,\Z)\to\SL(d,\Z/n\Z)\bigr).
\]
Let $I$ be a family of infinite subsets of prime numbers
such that $|i\cap j|<\infty$ for any $i,j\in I$ with $i\neq j$.
(We note that there exists such an uncountable family $I$.)
Associate each $i=\{p_1<p_2<\cdots\}\in I$ with the
decreasing sequence of finite index normal subgroups
$\G^{(i)}_n=\G_{i(n)}$ where $i(n)=p_1\cdots p_n$.
Then, the family $I$ is admissible.
\end{lem}

\begin{proof}
First, we note that $\G_m\cap\G_n=\G_{\gcd(m,n)}$.
By the celebrated results of Kazhdan for $d\geq 3$ (see \cite{bdhv})
and Selberg for $d=2$ (see \cite{lubotzky}) the group $\G$ has
the property $(\tau)$ with respect to the family $\{\G_n : n\in\N\}$.
We observe that the index $[ \G : \G^{(i)}_m\G^{(j)}_n ]$
is the cardinality of $\G$-orbits of
$(\G/\G^{(i)}_m)\times(\G/\G^{(j)}_n)$.
Since
\[
\SL(d,\Z/p_1\cdots p_l\Z)=\prod_{k=1}^l\SL(d,\Z/p_k\Z)
\]
for any mutually distinct primes $p_1,\ldots,p_l$,
one has a group isomorphism
\[
\SL(d,\Z/i(m)\Z)\times\SL(d,\Z/j(n)\Z)\cong\SL(d,\Z/k\Z)\times\SL(d,\Z/l\Z),
\]
where $k=\gcd(i(m),j(n))$ and $l=i(m)j(n)/\gcd(i(m),j(n))$.
Since
\[
(\G/\G^{(i)}_m)\times(\G/\G^{(j)}_n)\subset\SL(d,\Z/i(m)\Z)\times\SL(d,\Z/j(n)\Z)
\]
as a $\G$-set, one has
\[
[ \G : \G^{(i)}_m\G^{(j)}_n ]\le|\SL(d,\Z/k\Z)|\,[\SL(d,\Z/l\Z) : \G/\G_l].
\]
Therefore, the condition $\sup\{ [ \G : \G^{(i)}_m\G^{(j)}_n ] : m,n\in\N\}<\infty$
follows from the fact that $|i\cap j|<\infty$.
\end{proof}

For example, we can take $\G\le\SL(2,\Z)$ to be
$\ip{\left(\begin{smallmatrix} 1& 2\\ 0& 1\end{smallmatrix}\right),
\left(\begin{smallmatrix} 1& 0\\ 2&
1\end{smallmatrix}\right)}\cong\F_2$. By \cite{shalom}, one may
relax the assumption that ``$\G\le\SL(d,\Z)$ has finite index'' to
``$\G\le\SL(d,\Z)$ is co-amenable,'' so that one can take $\G$ to be
isomorphic to $\F_\infty$.

Let $\Sg=(\G_n)_{n=1}^\infty$ be a decreasing sequence of finite
index subgroups of a group $\G$.
We write $X_{\Sg}=\varprojlim\G/\G_n$ for the projective limit of
the finite probability space $\G/\G_n$ with uniform measures.
We note that $L^\infty(X_{\Sg})=(\bigcup\ell^\infty(\G/\G_n))''$,
where the inclusion
$\iota_n\colon\ell^\infty(\G/\G_n)\hookrightarrow\ell^\infty(\G/\G_{n+1})$
is given by $\iota_n(f)(g\G_{n+1})=f(g\G_n)$. There is a natural
action $\G\curvearrowright L^\infty(X_{\Sg})$ which is ergodic,
measure-preserving and profinite. (Any such action
arises in this way.) The action is essentially-free if and only if
\begin{equation}\label{eq:ess-free}
\forall g\in\G\setminus\{1\}\quad
|\{ x\in X_{\Sg} : gx=x\}|=\lim_n \frac{|\{ x\in\G/\G_n : gx=x\}|}{|\G/\G_n|}=0.
\end{equation}
This condition clearly holds if all $\G_n$ are normal and $\bigcap\G_n=\{1\}$.
We denote $A_{\Sg}=L^\infty(X_{\Sg})$ and
$A_{\Sg,n}=\ell^\infty(\G/\G_n)\subset A_{\Sg}$. Since
\[
L^2(A_{\Sg})\cong\C1\oplus\bigoplus_{n=1}^\infty
\bigl(L^2(A_{\Sg,n})\ominus L^2(A_{\Sg,n-1})\bigr)
\subset\C1\oplus\bigoplus_{n=1}^\infty\ell^2(\G/\G_n)^{o}
\]
as a $\G$-space, the action $\G\curvearrowright A_{\Sg}$ is
strongly ergodic if $\G$ has the property $(\tau)$ with respect to $\Sg$.

\begin{thm}\label{thm:uncountable}
Let $\G$ be a countable group which is not inner-amenable,
and $I$ be an uncountable admissible family of decreasing
sequences of finite index normal subgroups of $\G$.
Then, all $M_i=L(X_i)\rtimes\G$ are full factors of type
$\mathrm{II}_1$ and the set $\{ M_i : i\in I\}$ contains
uncountably many isomorphism classes of von Neumann algebras.
\end{thm}
\begin{proof}
That all $M_i$ are full follows from \cite{choda}.
Take a finite subset $F$ of $\G$ and $\e>0$
such that $(F,\e)$ witnesses spectral gap
for both non-inner-amenability and the property $(\tau)$
with respect to $\{\G^{(i)}_m\cap\G^{(j)}_n\}$.
We write $\lambda_i(g)$ for the unitary element in $M_i$
that implements the action of $g\in\G$.

We claim that if $i\neq j$,
then $(F,\e)$ witnesses essential spectral gap of
the unitary $\G$-representation
$\Ad(\lambda_i\otimes\lambda_j)$ on $L^2(M_i\vt M_j)$.
First, we deal with the $\Ad(\lambda_i\otimes\lambda_j)(\G)$-invariant
subspace
\begin{equation}\label{eq:AIAJ}
L^2(A_i\vt A_j)
\cong\C1\oplus\bigoplus_{n=1}^\infty\bigl(L^2(A_{i,n}\vt A_{j,n})
\ominus L^2(A_{i,n-1}\vt A_{j,n-1})\bigr).
\end{equation}
We note that the unitary $\G$-representation on
\[
L^2(A_{i,n}\vt A_{j,n})\cong\ell^2((\G/\G^{(i)}_n)\times(\G/\G^{(j)}_n))
\]
is contained in a multiple of $\ell^2(\G/(\G^{(i)}_n\cap\G^{(j)}_n))$.
Hence if we show that the subspace of
$\G$-invariant vectors in $L^2(A_i\vt A_j)$
is finite-dimensional, then we can conclude by the property $(\tau)$
that $(F,\e)$ witnesses essential spectral gap.
Suppose $\xi\in L^2(A_{i,n}\vt A_{j,n})$ is $\G$-invariant.
Since $\G^{(i)}_n$ acts trivially on $L^2(A_{i,n})$, the vector
$\xi$ is $\Ad(1\otimes\lambda_j)(\G^{(i)}_n)$-invariant.
The same thing is true for $j$.
It follows that $\xi$ is in the
$\G^{(i)}_n\G^{(j)}_n\times\G^{(i)}_n\G^{(j)}_n$-invariant
subspace, whose dimension is $[\G:\G^{(i)}_n\G^{(j)}_n]^2$.
Since this number stays bounded as $n$ tends to $\infty$, we are done.
Second, we deal with the
$\Ad(\lambda_i\otimes\lambda_j)(\G)$-invariant subspace
\begin{equation}\label{eq:AIMJ}
(L^2(M_i)\ominus L^2(A_i))\vt L^2(M_j)\cong
\ell^2(\G\setminus\{1\})\vt L^2(A_i)\vt L^2(M_j),
\end{equation}
where $\G$ acts on the right hand side Hilbert space
(which will be denoted by $\hh$) as
$\Ad(\lambda(g)\otimes\lambda_i(g)\otimes\lambda_j(g))$.
For every vector $\xi\in\hh$,
we write it as $(\xi_g)_{g\in\G\setminus\{1\}}$ with
$\xi_g\in L^2(A_i)\vt L^2(M_j)$ and
define $|\xi|\in\ell^2(\G\setminus\{1\})$ by $|\xi|(g)=\|\xi_g\|$.
It follows that
\begin{align*}
\Re\ip{\Ad(\lambda(g)\otimes\lambda_i(g)\otimes\lambda_j(g))\xi,\xi}
&= \Re\sum_{h\in\G\setminus\{1\}}
 \ip{\Ad(\lambda_i(g)\otimes\lambda_j(g))\xi_h,\xi_{ghg^{-1}}}\\
&\le \sum_{h\in\G\setminus\{1\}} \|\xi_h\|\|\xi_{ghg^{-1}}\|
=\ip{\Ad\lambda(g) |\xi|,|\xi|}
\end{align*}
for every $g\in\G$ and $\xi\in\hh$.
Since $(F,\e)$ witnesses spectral gap of the conjugation action
on $\ell^2(\G\setminus\{1\})$, it also
witnesses spectral gap of the $\G$-action on $\hh$.
Similarly, $(F,\e)$ witnesses spectral gap of
\begin{equation}\label{eq:MIAJ}
L^2(M_i)\vt (L^2(M_j)\ominus L^2(A_j)).
\end{equation}
Since the Hilbert spaces (\ref{eq:AIAJ}--\ref{eq:MIAJ}) cover
$L^2(M_i\vt M_j)$, we conclude that $(F,\e)$ witnesses
essential spectral gap of the $\G$-action
$\Ad(\lambda_i\otimes\lambda_j)$.
This argument is inspired by \cite{choda}.

We claim that for any $i\in I$ and any unitary element
$u(g)\in M_i$ with $\|\lambda_i(g)-u(g)\|_2<\e/4$,
the essential spectrum of the self-adjoint operator
\[
h_F=\frac{1}{2|F|}
\sum_{g\in F}\bigl(\Ad(\lambda_i(g)\otimes u(g))
+\Ad(\lambda_i(g^{-1})\otimes u(g^{-1}))\bigr)
\]
on $L^2(M_i\vt M_i)$ intersects with $[1-\e/2,1]$.
We fix $i\in I$ and define for every $n\in\N$
the projection $\chi_n\in M_i\vt M_i$ by
$\chi_n=\sum e_k\otimes e_k$, where $\{e_k\}$
is the set of non-zero minimal projections in
$A_{i,n}\cong\ell^\infty(\G/\G^{(i)}_n)$.
We normalize $\xi_n=[\G:\G^{(i)}_n]^{1/2}\chi_n$
so that $\|\xi_n\|_2=1$.
Then, it is not hard to see
\[
\Ad(\lambda_i(g)\otimes\lambda_i(g))\xi_n=\xi_n
\]
for all $g\in\G$, and
\[
\|(1\otimes a)\xi_n\|_2^2=\|a\|_2^2=\|\xi_n(1\otimes a)\|_2^2
\]
for all $a\in M_i$. It follows that
\begin{alignat*}{2}
\ip{h_F\xi_n,\xi_n}
&= \frac{1}{|F|}\sum_{g\in F}
 \Re\ip{\Ad(\lambda_i(g)\otimes u(g))\xi_n,\xi_n}\\
&\geq \frac{1}{|F|}\sum_{g\in F}(1-2\|\lambda_i(g)-u(g)\|_2)
& &>1-\e/2.
\end{alignat*}
Since $\xi_n\to0$ weakly as $n\to\infty$,
the claim follows (cf.\ \cite{ioana:cocycle}).

>From the above claims, we know that if $i\neq j$, then there is no
$*$-isomorphism $\theta$ from $M_i$ onto $M_j$ such that
$\|\theta(\lambda_i(g))-\lambda_j(g)\|_2<\e/4$ for all $g\in F$.
Now, if the isomorphism classes of $\{ M_i : i\in I \}$ were
countable, then there would be $M_0$ and an uncountable subfamily
$I_0\subset I$ such that $M_i\cong M_0$ for all $i\in I_0$. Take an
$*$-isomorphism $\theta_i\colon M_i\to M_0$ for every $i\in I_0$.
Since $M_0^{F}$ is separable in $\|\,\cdot\,\|_2$-norm, there has to
be $i,j\in I_0$ with $i\neq j$ such that
\[
\max_{g\in F}\|\theta_{i}(\lambda_{i}(g))
 -\theta_{j}(\lambda_{j}(g))\|_2<\e/4,
\]
in contradiction to the above.
\end{proof}

When combined with Lemma 5.3, Theorem \ref{thm:uncountable} shows in
particular that any arithmetic property (T) group has uncountably
many orbit inequivalent free ergodic profinite actions, thus
recovering a result in \cite{ioana:cocycle}. However,
\cite{ioana:cocycle} provides a ``concrete'' family (consequence of
a cocycle superrigidity result for profinite actions of Kazhdan
groups) rather than an ``existence'' result, as
\ref{thm:uncountable} does. But the consequence of
\ref{thm:uncountable}, 5.3 that's relevant here is the following:

\begin{cor}\label{cor:uncappfree}
For each $2\le r\le\infty$, there exist uncountably many
non-isomorphic approximate free group factors of rank $r$. In
particular, there exist uncountably many orbit inequivalent free
ergodic profinite actions of $\Bbb F_r$.
\end{cor}

\begin{rem}\label{rem:appfree} Note that if $2\le r\le\infty$ and
$\Sg=(\G_n)$ is a decreasing sequence of finite index subgroups of
the free group $\F_r$ satisfying the condition (\ref{eq:ess-free}),
then the associated free group factor of rank $r$ is the inductive
limit of $A_{\Sg,n}\rtimes\F_r\cong\B(\ell^2(\F_r/\G_n))\vt
L(\G_n)$, which is isomorphic to $L(\F_{1+(r-1)/[\G:\G_n]})$, by
Schreier's and Voiculescu's formulae (\cite{vdn}). Since
$1+(r-1)/[\G:\G_n] \rightarrow 1$, this justifies the notation
$L(\F_{\oneplus}^{r,\Sg})$ for the approximate free group factor
$L^\infty(X_{\Sg})\rtimes\F_r$. The factors $L(\F_{\oneplus}^*)$ can
be viewed as complementing the one parameter family of free group
factors $L(\F_{1+t}), 0 < t \leq \infty$, in
\cite{dykema,radulescu}.

As mentioned in Section 4, all $L(\F_{\oneplus}^{r,\Sg})$ have
Haagerup's compact approximation property (by \cite{haagerup:map}),
the complete metric approximation property (by \ref{Lambdacb}) and
unique Cartan subalgebra, up to unitary conjugacy (by
\ref{cor:uniquecartan}). Also, by \cite{ozawa:kurosh}, the commutant
of any hyperfinite subfactor of $L(\F_{\oneplus}^{r,\Sg})$ must be
an amenable von Neumann algebra, in particular
$L(\F_{\oneplus}^{r,\Sg})$ is prime, i.e.\ it cannot be written as a
tensor product of two II$_1$ factors. By \cite{popa:betti}, since
the factors $L(\F_{\oneplus}^{r,\Sg})$ have Haagerup property they
cannot contain factors $M$ which have a diffuse subalgebra with the
relative property (T). In particular, the $\mathcal H\mathcal
T$-factors considered in \cite{popa:betti}) cannot be embedded into
approximate free group factors. Same for the factors arising from
Bernoulli actions of ``w-rigid'' groups in \cite{popa:bernoulli}.

Corollary \ref{cor:uniquecartan} combined with \cite{gaboriau} shows
that approximate free group factors of different rank are
non-isomorphic, $L(\F_{\oneplus}^{r,\Sg})\not\simeq
L(\F_{\oneplus}^{s,\Sg})$, $\forall 2\leq r \neq s \leq\infty$, and
have trivial Murray-von Neumann fundamental group \cite{mvn} when
the rank is finite, $\mathcal{F}(L(\F_{\oneplus}^{s,\Sg}))=\{1\}$,
$\forall 2\leq r < \infty$. (Recall from \cite{mvn} that if $M$ is a
II$_1$ factor then its fundamental group is defined by
$\mathcal{F}(M)=\{t> 0 \mid M^t \simeq M\}$.) The first examples of
factors with trivial fundamental group were constructed in
\cite{popa:betti}, were it is shown that
$\mathcal{F}(L^\infty({\mathbb T}^2)\rtimes \F_r)=\{1\}$, for any
finite $r \geq 2$, the action of $\F_r$ on ${\mathbb T}^2$ being
inherited from the natural action $\SL(2,\Z) \curvearrowright
{\mathbb T}^2=\hat{\Z^2}$, for some embedding
$\F_r\subset\SL(2,\Z)$.

One can show that amplifications of approximate free group factors
are related by the formula $L(\F_{\oneplus}^{r,\Sg})^{t} =
L(\F_{\oneplus}^{r',{\Sg}'})$, with $r'=t^{-1}(r-1)+1$, whenever
$t^{-1}$ is an integer dividing the index of some
$[\Gamma:\Gamma_n]$ in the decreasing sequence of groups
$\Sg=(\Gamma_n)$, with ${\Sg}'$ appropriately derived from $\Sg$. It
is not clear however if this is still the case for other values of
$t$ for which $t^{-1}(r-1)+1$ is still an integer.

Finally, note that $L(\F_{\oneplus}^{r,\Sg})$ is non-$\Gamma$ if and
only if the action $\Gamma \curvearrowright X_{\Sg}$ has spectral
gap. Indeed, since the acting group is $\F_r$, any asymptotically
central sequence in
$L(\F_{\oneplus}^{r,\Sg})=L^\infty(X_{\Sg})\rtimes \F_r$ must lie in
$L^\infty(X_{\Sg})$, so $L(\F_{\oneplus}^{r,\Sg})$ is non-$\Gamma$
if and only if $\F_r \curvearrowright X_{\Sg}$ is strongly ergodic,
which by \cite{abert-elek} is equivalent to $\F_r \curvearrowright
X_{\Sg}$ having spectral gap. For each $2 \leq r \leq \infty$, one
can easily produce sequences of subgroups $\Sg=(\Gamma_n)$ such that
$\F_r \curvearrowright X_{\Sg}$ does not have spectral gap, thus
giving factors $L(\F_{\oneplus}^{r,\Sg})$ with property $\Gamma$. On
the other hand, as mentioned before, if $\F_r$ is embedded with
finite index in $\SL(2,\Z)$ (or merely embedded ``co-amenably,'' see
\cite{shalom}) and $\Sg=(\Gamma_n)$ is given by congruence
subgroups, then $\F_r \curvearrowright X_{\Sg}$ has spectral gap by
Selberg's theorem. Thus, the corresponding approximate free group
factors $L(\F_{\oneplus}^{r,\Sg})$ are non-$\Gamma$. By Corollary
\ref{cor:uncappfree} and its proof, there are uncountably many
non-isomorphic such factors $L(\F_{\oneplus}^{r,\Sg})$ for each
$2\leq r \leq \infty$. It is an open problem on whether
there exist solid factors within this class.
\end{rem}

\end{document}